\documentclass[11pt,a4paper]{amsart}
\usepackage[margin=3cm]{geometry}

\usepackage[latin1]{inputenc}
\usepackage{amsmath,amsfonts,amssymb,amsthm}
\usepackage{graphicx}

\usepackage{enumerate,hyperref} 
\usepackage[version=3]{mhchem}
\usepackage[all]{xy}
\usepackage{adjustbox}
\usepackage{arydshln}
\usepackage{marginnote}
\usepackage{multirow}
\usepackage{resizegather}
\usepackage[dvips,pdftex]{color}
\usepackage{colortbl}
\usepackage{pgf,tikz}
\usepackage{url}
\usepackage{longtable}
\usepackage[mathlines]{lineno}
 
\usepackage[sort&compress,comma,square,numbers]{natbib}

\theoremstyle{plain}
\newtheorem{thm}{Theorem}[section]

\newtheorem{lem}[thm]{Lemma}

\theoremstyle{definition}

\newtheorem{example}[thm]{Example}

\newcommand{\cS}{\mathcal{S}}
\newcommand{\R}{\mathbb{R}}
\newcommand{\N}{{\mathbb N}}
\newcommand{\Z}{{\mathbb Z}}

\newcommand{\E}{{\mathbb E}}
\renewcommand{\k}{\kappa}

 \DeclareMathOperator{\diag}{diag}

\begin{document}
\author[E. Feliu, A. Sadeghimanesh]{Elisenda Feliu$^1$, AmirHosein Sadeghimanesh$^{2}$}
\title[Kac-Rice formulas and parameter regions]{Kac-Rice formulas and the number of solutions of parametrized systems of polynomial equations}

\footnotetext[1]{Department of Mathematical Sciences, University of Copenhagen. \emph{Address: } Universitetsparken 5, 2100 Copenhagen, Denmark.  \emph{Email: } efeliu@math.ku.dk. }
\footnotetext[2]{Centre for Computational Science \& Mathematical Modelling (CSM), Coventry University. \emph{Address: } The Futures Institute, Unit 10, Coventry Innovation Village, Cheetah Road, CV1 2TL Coventry, UK   \emph{Email: }  ad6397@coventry.ac.uk}

	\maketitle
	\begin{abstract}
Kac-Rice formulas express the expected number of elements a fiber of a  random field has in terms of a multivariate integral. We consider here parametrized systems of polynomial equations that are linear in enough parameters, and provide a Kac-Rice formula for the expected number of solutions of the system when the parameters follow continuous distributions. Combined with Monte Carlo integration, we apply the formula to partition the parameter region according to the number of solutions or find a region in parameter space where the system has the maximal number of solutions. The motivation stems from the study of steady states of chemical reaction networks and gives new tools for the open problem of identifying the parameter region where the network has at least two positive steady states. We illustrate with numerous examples that our approach successfully handles a larger number of parameters than  exact methods.


\end{abstract}

	\section*{Introduction}
Systems of parametrized polynomial equations arise naturally in applications, and in particular in relation to steady states of polynomial ordinary differential equations (ODEs). We address here the problem of describing the function mapping a parameter vector to the number of solutions of the  system specialized to the parameter vector.  
That is, given a parametrized system of $n$ polynomial equations in $n$ variables
\begin{equation}\label{eq:mysystem} f_\k(t)=0,\qquad t\in A,\quad \k\in B,
\end{equation}
with $A\subseteq \R^n$ and $B\subseteq \R^m$, we want to partition the parameter space $B$ into regions where the number of solutions to the system in $A$ is $0,1,2,\dots,+\infty$. 

The motivation stems from the study of the steady states of (bio)chemical reaction networks, where $A=\R^n_{>0}$ (the positive orthant), and $B$ typically is   $\R^m_{>0}$. In this setting, it is in particular of interest to understand for what parameter values the system describing the steady states of the network has at least two positive solutions (see Subsection~\ref{sec:Reaction_networks_and_multistationarity}). This property is termed \emph{multistationarity}, and implies that the network, corresponding to a chemical system or mechanism in the cell, for instance, can potentially rest in two different states under exactly the same conditions.  This property has received substantial attention in the context of systems and synthetic biology, for the implications in cell decision making \cite{BiochemicalSwitchInMammalianCells,StemCellSwitch,BistabilityHumanCancer}. In this context, only stable steady states are relevant, which leads to the concept of \emph{bistability}, referring to the existence of two stable positive steady states. This usually implies that the network has at least three positive steady states, two of which are stable and one is unstable. 
There exist numerous approaches to determine whether multistationarity exists for some choice of parameter values, e.g.   \cite{FeinbergExistenceAndUniquness,Feinberg1988,crnttoolbox,Dickenstein-Toric,Feliu-Sign,conradi-switch,feliu-bioinfo,control,AtomsJoshi,Feliu-Simplifying,Multistationarity}. However, finding the parameter regions where the system displays multistationarity is a much harder question. Only very recently, approaches to partially understand the region of multistationarity have  been proposed, e.g.~\cite{maya-bistab,FeliuPlos,dickenstein:regions}. 
Furthermore, most methods to identify multistationarity return a parameter value in $\R^m_{>0}$ for which multistationarity occurs, but cannot be adapted to determine whether multistationarity also occurs for parameters in a given subset of $\R^m_{>0}$. This has the consequence that often, returned parameter values do not fall into biochemically relevant regions.

Theoretically, in order to partition the parameter space according to the number of solutions of the system, or to identify the parameter region of multistationarity, one might employ quantifier elimination or Cylindrical Algebraic Decomposition (CAD)  \cite{basu2007algorithms}. However, these approaches have a high computational cost and become unfeasible already for small systems with three or four parameters and three or four variables.

In this work we explore the use of Kac-Rice formulas to study the parameter space.  Kac introduced in \cite{KacOriginalPaper} a formula  to compute the expected number of real roots of a univariate polynomial with random coefficients. At about the same time, Rice introduced a similar formula for the number of crossings of ergodic stationary processes~\cite{Rice:first,Rainal:OriginRice}. 
These formulas became known as Kac-Rice formulas, and were extended in several directions, e.g. \cite{Edelman_Kostlan,LevelSetsAndExtremaOfRandomProcessesAndFields,ylvisaker1965}, see also \cite{OnTheKacRiceFormula}. This led to a metaformula for the expected number of elements of a fiber of 
a random field on a manifold under some conditions \cite{RandomFieldsAndGeometry}. In general, a Kac-Rice formula expresses the expected number  by means of a multivariate integral, and has found applications in many areas such as regression \cite{kac-rice-inference},  the theory of random matrices \cite{auffinger},  number theory \cite{evans} or enumerative geometry \cite{basu-enumerative}, to name a few.

In the first part of this work, we derive a Kac-Rice formula suited to polynomial systems with ``suficient'' linearity in the parameters (Theorem~\ref{thm:Kac-Rice}), which accommodates systems arising from reaction networks.  The formula expresses, in terms of a multivariate integral, the expected number of solutions of system \eqref{eq:mysystem}, when the parameters $\k_1,\dots,\k_m$ are independent random variables with continuous distributions.
We provide a direct proof of the formula that combines usual arguments in the derivation of Kac-Rice formulas \cite{LevelSetsAndExtremaOfRandomProcessesAndFields}. 
We proceed to discuss how the Kac-Rice integral can be computed by employing numerical integration, mainly uisng Monte Carlo methods. 

In the second part of this work, we provide numerous examples to illustrate how the Kac-Rice formula can be used in practice to address the following problems:
\begin{enumerate}[(i)]
	\item Provide a grid partition of the parameter region $B$ according to the number of solutions of system~\eqref{eq:mysystem}.
\item Find a parameter point or region for which the system has the maximal number of solutions, or at least $M$ solutions (for some number $M$).
\end{enumerate}

These questions are addressed by endowing the parameters with the uniform distribution in a box (product of intervals). Then the Kac-Rice formula gives the average number of solutions the system has when the parameters belong to the box. By making the boxes small, we can approximately partition  the parameter region according to the number of solutions. 
We show that our approach can handle systems with over $10$ parameters, where exact methods fail due to computational power. We also illustrate how parallelisation of our computations enables the study of complex systems.

The organization of the paper is as follows. In Section \ref{sec:Kac-Rice-theory} the motivational setting of reaction networks is introduced, the statement of the Kac-Rice theorem is given, and  we discuss Monte Carlo integration.  Section~\ref{sec:examples} 
devises the strategy to use the Kac-Rice formula to study parameter regions, and illustrates it with numerous examples. Finally, Section~\ref{sec:proofKacRice} contains the proof of Theorem~\ref{thm:Kac-Rice}. 
Computational files can be accessed in the Github repository \cite{MCKR_project_repository}. A Julia implementation of the  methods in this manuscript for user-created input files is available at \cite{MCKR_program}. 

	\medskip
\noindent
	\textbf{Notation. }
 $\R_{\geq 0}$ and $\R_{>0}$ refer to the non-negative and positive real numbers respectively. A box $B\subseteq\R^n$ is a Cartesian product $\prod_{i=1}^nB_i$ of (possibly unbounded) intervals of the real line. The intervals are allowed to be (half-)open or closed.
For a set $B\subseteq \R^n$, we let $\chi_{B}\big(y\big)$ denote the indicator function being $1$ if $y\in B$ and $0$ otherwise. 
For a finite interval $B\subseteq \R$,  $U(B)$ denotes the uniform distribution on $B$. We simply write $U(a,b)$, if $a,b$ are the extremes of the interval $B$. 
We let $\bar{N}_{(a,b)}(\mu,\sigma)$ denote the truncated normal distribution in the interval $(a,b)$, derived from a normal distribution with mean $\mu$ and standard deviation $\sigma$. In this work, we take $\mu=\tfrac{a+b}{2}$, so $\mu$ is also the mean of the truncated normal distribution.

\section{Expected number of solutions using Kac-Rice formulas}\label{sec:Kac-Rice-theory}

\subsection{Motivation: Reaction networks and multistationarity}\label{sec:Reaction_networks_and_multistationarity}
In this section we introduce the polynomial system of interest in the study of steady states of reaction networks. 
A \emph{reaction network} on a set $\cS=\{X_1,\dots,X_n\}$ (\emph{species set}) is a collection of \emph{reactions} between  linear combinations of  species:
\begin{equation}\label{eq:reaction}
  \sum_{i=1}^n a_{ij} X_i \rightarrow \sum_{i=1}^n b_{ij} X_i, \qquad j=1,\dots,r,
\end{equation}
with $a_{ij},b_{ij}\in \Z_{\geq 0}$.
Let $x_i(t)$ denote the concentration of $X_i$ at time $t$ and $x(t)=(x_1(t),\dots,x_n(t))$. 
Under the so-called \emph{mass-action assumption} \cite{feinberg-book,gunawardena-notes}, the evolution of the concentrations of the species in time is modeled by means of a polynomial system of autonomous ODEs in $\R^n_{\geq 0}$ of the form:
\begin{equation}\label{eq:ode} 
\frac{dx(t)}{dt} =F_{k}(x(t)), \qquad \textrm{where } F_{k,i}(x)= \sum_{j=1}^r (b_{ij}-a_{ij}) k_{j}\,  x_1^{a_{1j}}\cdots x_n^{a_{nj}},\quad i=1,\dots,n.
\end{equation}
Here $k_j>0$ are called \emph{reaction rate constants}, and $0^0=1$ by convention. Typically, $k_j$ are considered labels of the reactions, and by default their subindex indicates the order of the set of reactions.
By letting $N\in \Z^{n\times r}$ be the matrix with entries $b_{ij}-a_{ij}$ for $i=1,\dots,n$, $j=1,\dots,r$,  any vector $\omega$ in the left kernel of $N$ gives rise to a linear first integral, as $\omega\cdot   \frac{dx(t)}{dt}=0$.  Hence there are invariant linear subspaces with equations 
\[W x = T,\qquad T\in \R^d,\]
for any matrix $W\in \R^{d\times n}$ whose rows form a basis of $\ker(N^t)$.
These equations are called \emph{conservation laws}, and $T$ a vector of \emph{total amounts}.

\medskip
The \emph{steady states} of the ODE system \eqref{eq:ode} in the invariant linear subspace with total amount $T$ are the non-negative solutions to the system $F_{k}(x)=0$, $W x - T=0$.
As the conservation laws describe linear relations among the entries of $F_k$,   $d$ entries of $F_k$  are redundant (linearly dependent of the rest) and can be removed.  Let $\widetilde{F}_{k}(x)$ be a function with $n-d$ entries obtained in this way. Then the system of interest is square with $n$ variables and $n$ equations:
\begin{equation}\label{eq:steadystates}
\widetilde{F}_{k}(x)=0, \qquad W x - T=0.
\end{equation}

The network is said to be \emph{multistationary} if there exist $k=(k_1,\dots,k_r)\in \R^r_{>0}$ and $T\in \R^d$ such that system \eqref{eq:steadystates} admits at least two positive solutions. 
Our ultimate goal is to understand how the number of positive solutions to \eqref{eq:steadystates} depends on $k\in \R^r_{>0}$ and  $T\in \R^d$.
This implies understanding the following map:
\begin{eqnarray}
\R^{r}_{>0} \times \R^d & \rightarrow & \N \cup \{+\infty\}  \label{eq:map_parameters} \\
\k:=(k,T) & \mapsto &  \#  \{ x\in \R^n_{>0} \mid x \textrm{ is a solution to \eqref{eq:steadystates}}   \}. \nonumber
\end{eqnarray}
The image of this map partitions the parameter space $ \R^{r}_{>0} \times \R^d$. 

\subsection{The Kac-Rice formula}	
We give here a Kac-Rice formula on the expected number of solutions of a polynomial system, which applies to polynomial systems with ``sufficient'' linearity in the parameters. This will later be applied to understand the map \eqref{eq:map_parameters}. 
In the following, measurability is with respect to the Borel $\sigma$-algebra on $\R^n$, 
and integrals are considered with respect to the Lebesgue measure on $\R^n$.  
 
We consider functions of $n$ polynomials  in $n$ variables and $m$ parameters
\begin{equation*}
f_\k(t)=(f_{\k,1}(t),\dots,f_{\k,n}(t)),\qquad t=(t_1,\dots,t_n)\in \R^n,\quad  \k=(\k_1,\dots,\k_m)\in \R^m,
\end{equation*}
with $m\geq n$, and such that  the coefficients of the entries of $f_\k(t)$ are polynomials in $\k$.  
In the motivating scenario from Subsection~\ref{sec:Reaction_networks_and_multistationarity}, the polynomial map is given by the left-hand side of \eqref{eq:steadystates} and $\k=(k,T)$ such that $m=r+d$.

We assume that the parameters $\k_1,\dots,\k_m$ are independent random variables with continuous distribution and  density $\rho_i$ in an interval $B_i$, for $i=1,\dots,m$. If $t$ is such that the image of $f_\k(t)$ has non-zero   measure in $\R^n$, then the values of $f_\k(t)$ for varying $\k$ 
define a random  variable taking values in $\R^n$ with a continuous distribution induced by the densities $\rho_i$.  
As  $f_\k(t)$ is polynomial in $\k$, this image has non-zero measure if and only if it is Zariski dense in $\R^n$, or equivalently, the Jacobian of the polynomial map $f_\k(t)\colon \R^m \rightarrow \R^n$ with variable $\k$ has maximal rank $n$.
In particular,   the image of $f_\k(t)$ in $\k$  neither is constant nor lies in a proper algebraic variety of $\R^n$.

\medskip
For a subset $A\subseteq\R^n$, consider the discrete random variable 
with state space $\Z_{\geq 0} \cup \{+\infty\}$ that counts the number of zeroes of $f_\k$ in $A$, and let $\mathbb{E}\big( \#(f_\k^{-1}(0)\cap A)\big)$ be its expected value. 
We let $J_{f_\k}(t) = \big( \frac{\partial f_{\k,i}(t)}{\partial t_j} \big)_{i,j} \in \R^{n\times n}$ be the Jacobian matrix of $f_\k(t)$ with respect to $t$.

The following theorem gives a Kac-Rice formula for $\mathbb{E}\big( \#(f_\k^{-1}(0)\cap A)\big)$ for polynomial functions of a certain form, in line with the Kac-Rice metatheorem from \cite[Ch 11]{RandomFieldsAndGeometry}. 
The proof is given in Section~\ref{sec:proofKacRice}.

\begin{thm}[Kac-Rice formula]\label{thm:Kac-Rice}
Let $A=A_1\times\cdots\times A_n\subseteq \R^n$ be a box.  
Let $f_\k\colon A\rightarrow\mathbb{R}^n$ be a polynomial map whose coefficients are polynomials in $\k=(\k_1,\dots,\k_m)$ with $m\geq n$. 
Assume that each parameter $\k_i$ follows a continuous distribution with support on an interval $B_i$ and density $\rho_i$, such that $\k_1,\dots,\k_m$ are independently distributed. Assume $\rho_i$ is a continuous function except maybe for a finite number of points of $B_i$, for $i=1,\dots,n$.

Define
$\widetilde{B}=B_{n+1}\times\cdots\times B_m$ and let $\bar{\k}=(\k_{n+1},\ldots,\k_m)$. 
For each $i=1,\dots,n$, assume that there exist polynomials $h_i(\bar{\k},t)$ and $q_i(\bar{\k},t)$ in $\bar{\k},t$,  such that 
\begin{equation}\label{eq:decompose:f}
f_{\k,i}(t)= h_i(\bar{\k},t)\k_i+q_i(\bar{\k},t).
\end{equation}
For $(\bar{\k},t)\in \widetilde{B} \times A$ define  
	\[   g_{\bar{\k},i}(t) :=\frac{-q_i(\bar{\k},t)}{ h_i(\bar{\k},t)},\quad i=1,\dots,n, \qquad g_{\bar{\k}}(t):=( 
	g_{\bar{\k},1}(t) , \dots, g_{\bar{\k},n}(t)), \]
	and
	\[\bar{\rho}(\bar{\k},t):=
	\left(\prod_{i=1}^n \rho_i(g_{\bar{\k},i}(t)) \right)\left(\prod_{i=n+1}^m \rho_i(\k_i)\right)\quad \textrm{if }\bar{\k}\in \widetilde{B},\qquad \textrm{and}\quad\bar{\rho}(\bar{\k},t):=0\quad \textrm{otherwise}. \]

Assume that  
\begin{enumerate}[(i)]
\item  $h_i(\bar{\k},t)$ does not vanish in $\widetilde{B}\times A$ for $i=1,\dots,n$.
\item For $\bar{\k}$ outside a Zariski closed set (relative to $\widetilde{B}$) of  measure zero $\widetilde{P}\subseteq \widetilde{B}$, 
the numerator of $\det(J_{g_{\bar{\k}}}(t))$ is a non-zero polynomial in $t$; equivalently $\det(J_{f_\k}(t))_{|(\k_1,\dots,\k_n)=g_{\bar\k}(t)}\neq 0$ as a rational function in $t$, see \eqref{eq:Jacfg}.
\end{enumerate}

	Then for all $t\in A$, the image of $f_\k(t)$ has positive measure in $\R^n$ and 
			\begin{align}\label{eq:kac-rice-formula}
 	\mathbb{E}(\# (f_\k^{-1}(0) \cap A)) & = \int_A \int_{\widetilde{B}} \  \big| \det\big(J_{g_{\bar{\k}}}(t)\big)\big|\, \bar{\rho}(\bar{\k},t)   \, d\k_{n+1}\dots d\k_m
   \, dt.
	\end{align}

	\end{thm}
Equality \eqref{eq:kac-rice-formula}  is  called the \emph{Kac-Rice formula}, and the integral on the right-hand side of the equality is called the \emph{Kac-Rice integral}. 	
 
The proof of Theorem~\ref{thm:Kac-Rice} is given in Section~\ref{sec:proofKacRice}. 
We typically consider  uniform or normal distributions on the parameters, hence the density functions are continuous outside a finite number of points. Note that Theorem~\ref{thm:Kac-Rice}(ii) implies that $q_i(\bar{\k},t)\neq 0$ as a polynomial in $\bar{\k},t$ for $i=1,\dots,n$. An easy computation shows that  for $\k\in B_1\times \dots \times B_m$, we have
\begin{equation}\label{eq:Jacfg}
\det\big(J_{g_{\bar{\k}}}(t)\big)= \tfrac{(-1)^{n}}
{\prod_{i=1}^n h_i(\bar{\k},t) } \det(J_{f_\k}(t))_{|(\k_1,\dots,\k_n)=g_{\bar\k}(t)}. 
\end{equation}
Using this, one can show that \eqref{eq:kac-rice-formula} agrees with the usual expression of Kac-Rice formulas, 
\[ \int_A \mathbb{E}\big(|\det(J_{f_\k}(t))|\mid f_\k(t)=0\big)p_t(0) dt,\] where $p_t(0)$ is the density of the random variable $f_\k$ at $0$. 
However, derivation of the formula in the form \eqref{eq:kac-rice-formula} is more straightforward and avoids considering random variables conditioned on a  measure zero set.

\medskip
Back to the motivating scenario, the next theorem guarantees that the Kac-Rice formula can be applied to study  systems arising from reaction networks as in \eqref{eq:steadystates}.

\begin{thm}\label{Theorem:Kac-Rice_in_CRN}
System \eqref{eq:steadystates} is equivalent to a system that admits a decomposition of the form \eqref{eq:decompose:f} satisfying assumption (i) of Theorem~\ref{thm:Kac-Rice} for $A\subseteq \R^n_{>0}$ and any choice of intervals $B_i\subseteq \R$ for $i=1,\dots,m$. 
\end{thm}
\begin{proof}
Let $\k=(k,T)$. 
The $i$-th equation in $Wx-T=0$ (indexed $n-d+i$ in  \eqref{eq:steadystates}) decomposes as in \eqref{eq:decompose:f} with parameter $T_i$ and $h_{n-d+i}(\k,t)=1$. 

Consider now the entries of $\widetilde{F}_k$ and let $\widetilde{F}_k = \widetilde{N} \diag(k) x^Y$, where $Y=(a_{ij})\in \R^{n\times r}$ is the matrix of coefficients of the reactants, see \eqref{eq:reaction}. By construction, $\widetilde{N}$ is any choice of $n-d$ linearly independent rows of $N$. 
Hence there exist column indices $i_1,\dots,i_{n-d}$ such that the submatrix $\widetilde{N}'$ of $\widetilde{N}$ given by these columns has full rank $n-d$. For simplicity assume $i_j=j$, and write $\widetilde{N}= (\widetilde{N}' | \widetilde{N}'')$.  
Consider the function 
\[G_k(x)=  (\widetilde{N}')^{-1} \widetilde{F}_k(x) = \big(\, {\rm Id}_{n-d} \, | \, (\widetilde{N}')^{-1}\widetilde{N}'' \, \big)  \diag(k) x^Y . \]
The solutions to $\widetilde{F}_k(x)=0$  and to $G_k(x)=0$ agree. 
Furthermore, $G_k(x)$ admits a decomposition as in \eqref{eq:decompose:f} with $\bar\k=(k_{1},\dots,k_{n-d})$, 
$h_i(\bar\k,x)= x_1^{a_{1i}} \cdots x_n^{a_{ni}}$ and $q_i(\bar\k,x)$ the $i$-th row of $(\widetilde{N}')^{-1}\widetilde{N}'' \diag(\bar\k) x^{\bar{Y}}$, with $\bar{Y}$ consisting of the last $r-(n-d)$ columns of $Y$. 
Clearly, $h_i(\bar\k,x)$ does not vanish in $\R^n_{>0}$. 
\end{proof}

We illustrate  Theorem~\ref{thm:Kac-Rice} with a couple of simple examples, before we turn to computing the Kac-Rice integral.

\begin{example}\label{Example:simplest}
Perhaps the simplest non-trivial example to consider is the linear polynomial $f_\k(t)=\k_2 t-\k_1$, which has one positive root if $\k_1\k_2>0$. Assume $\k_1,\k_2$ follow a uniform distribution in $[0,1]$ and  that $t$ is positive (that is, $A=\R_{>0}$). Note that in this case $\mathbb{E}\big(\#(f_{\kappa}^{-1}(0)\cap \mathbb{R}_{>0})\big)=1$, since the system has one positive solution for almost all $\kappa\in[0,1]^2$. We apply Theorem~\ref{thm:Kac-Rice}, with $g_{\k_2,1}(t)=\k_2 t$. Assumption (ii) holds with $\widetilde{P}=\{\k\in [0,1]^2 \mid \k_2 = 0\}$. We obtain 
	\begin{align*}
	\mathbb{E}\big( \#(f_\k^{-1}(0)\cap \R_{>0})\big) & = \int_0^{+\infty}\int_0^1  \k_2\, \chi_{[0,1]}(\k_2 t)d\k_2 \, dt
	= \int_0^{+\infty}\int_0^{\min(1,\frac{1}{t})} \k_2 d\k_2 \, dt  \\ & =
	 \int_0^{+\infty}  \frac{\min(1,\tfrac{1}{t})^2}{2}  dt = \int_0^1\tfrac{1}{2} dt+\int_1^{+\infty}\tfrac{1}{2t^2} dt=1,	\end{align*}
which is the correct value.
\end{example}

\begin{example}\label{Example:simple2v3p}
Consider this simple system of polynomial equations $f_\k(t)=0$
\[ \k_1 - \k_3 t_1 =0,\qquad \k_2-\k_3 t_1t_2=0.\]
We use Theorem~\ref{thm:Kac-Rice}
to determine the average number of solutions $(t_1,t_2)\in [0,1]^2$ in the parameter box $[0,1]^3$. 
With the notation of Theorem~\ref{thm:Kac-Rice}, we have $\bar{\k}=(\k_3)$, $g_{\k_3,1}(t)= \k_3t_1$ and
$g_{\k_3,2}(t)= \k_3t_1t_2$, with $h_1=h_2=1$. We consider each $\k_1,\k_2,\k_3$ uniformly distributed in $[0,1]$.  We have 
$|\det(J_{g_{\bar\k}}(t))| = \k_3^2t_1$. Hence assumption Theorem~\ref{thm:Kac-Rice}(ii) holds  with $\widetilde{P}=\{ \k \in [0,1]^3 \mid \k_3= 0\}$. 
This leads to the following:
\begin{align*}
\mathbb{E}\big(\#(f_\k^{-1}(0)\cap [0,1]^2)\big) &= \int_0^1\int_0^1 \int_0^1 \k_3^2t_1 \rho_1(\k_3t_1)\rho_2(\k_3t_1t_2) \rho_3(\k_3)d\k_3dt_1dt_2 \\ &= \int_0^1\int_0^1 \int_0^1 \k_3^2t_1 \, d\k_3dt_1dt_2 = \tfrac{1}{6}, 
\end{align*}
where we have used $\k_3t_1<1$ and $\k_3t_1t_2<1$.

Note that the solution to the system is $t_1=\tfrac{\k_1}{\k_3}, t_2=\tfrac{\k_2}{\k_1}$. This solution belongs to $[0,1]$ if and only if $\k_2<\k_1<\k_3$. The volume of this region within the cube $[0,1]^3$ is precisely $\tfrac{1}{6}$, in accordance with the result given by the Kac-Rice integral. 

\medskip
If the second equation is replaced with $\k_2t_2-\k_3 t_1t_2=0$, then  
$h_2(\k_3,t)=t_2$ vanishes in $A$ and hence Theorem~\ref{thm:Kac-Rice} does not apply. However, after factoring this equation as $t_2(\k_2 - \k_3t_1)$, the set of solutions of the original system is the union of the solution sets of two systems, arising from each factor, and for each of these systems Theorem~\ref{thm:Kac-Rice} applies.
\end{example}

\subsection{Monte Carlo integration for the Kac-Rice formula}\label{sec:Numerical_integration}
Although in some cases, such as in Examples~\ref{Example:simplest} and \ref{Example:simple2v3p},  the exact value of the Kac-Rice integral can be found, this is  typically not the case and one needs to rely on numerical integration.  
To this end,  we use Monte Carlo  integration with importance sampling.

\medskip
\paragraph{\bf Monte Carlo integration. }
We give the main ingredients of Monte Carlo integration relevant to this work  (see \cite{MonteCarloBook} for details). 
We consider an integral on a region $M\subseteq\R^{u}$ of the form
\begin{align}\label{eq:I}
I &= \int_{M}f(x_1,\dots,x_{u})dx_1\dots dx_{u}.
\end{align}
Given any probability distribution $P$   with non-zero probability density function $p(x)=p(x_1,\dots,x_{u})$ on $M$, it holds 
\begin{align*}
I &=\int_M \tfrac{f(x)}{p(x)}\, p(x)\, dx=\mathbb{E}\left(\tfrac{f(x)}{p(x)}\right).
\end{align*}
That is, the integral is expressed as the expected value of the function $Q(x):=\tfrac{f(x)}{p(x)}$ with respect to the chosen probability distribution. 
By the Law of Large Numbers, for large $N\in \N$, the integral $I$ can be approximated by the average value of $Q$ evaluated at randomly sampled points  $x^{(1)},\dots,x^{(N)}$ from $P$, that is, by
\begin{align}\label{Eq:Integral_Monte_Carlo}
\widehat{I}_N &= \tfrac{1}{N}\sum_{i=1}^N Q(x^{(i)}).
\end{align}
Furthermore, an estimate of the standard error of the approximation  is
\begin{align}\label{Eq:Standard_error_Monte_Carlo}
\hat{e}_N &   =  \sqrt{\tfrac{\sum_{i=1}^N\big(Q(x^{(i)}) -\widehat{I}_N\big)^2}{N(N-1)}}
=  \sqrt{\frac{\frac{1}{N}\sum_{i=1}^NQ(x^{(i)})^2-\widehat{I}_N^2}{N-1}},
\end{align}
where the second equality is well known (and easy to derive, see \cite{MonteCarloBook}).

We apply the approximation in \eqref{Eq:Integral_Monte_Carlo} to the Kac-Rice integral $I$ in \eqref{eq:kac-rice-formula} of Theorem~\ref{thm:Kac-Rice}.
To this end, we need  to choose a probability distribution on the domain $M=A\times \widetilde{B}  \subseteq \R^n\times \R^{m-n}$. 
In the applications in Section~\ref{sec:examples},  the box $A$ is bounded (or there exists a bounded set $A'\subseteq A$ containing $f_\k^{-1}(0)$ for all $\k\in B_1\times \dots \times B_m$), and we simply sample 
$t$ using the uniform distribution on $A$.
Let $\mu(t)$ denote the density of the  chosen distribution for $t$. 

For the integral over $\widetilde{B}$ in the $m-n$ parameters $\k_{n+1},\dots,\k_m$, we simply use the original density function $\rho_{n+1}\times \cdots \times \rho_m$. This choice makes the expression of the corresponding sums  in \eqref{Eq:Integral_Monte_Carlo}  simpler,  thereby   increasing the  computational speed. 
Specifically, with these choices, the function $Q=f/p$ used in \eqref{Eq:Integral_Monte_Carlo} becomes
	\begin{linenomath}
\begin{align*}
Q(t,\bar\k) &= \big| \det\big(J_{g_{\bar{\k}}}(t)\big) \big|\, \frac{\bar{\rho}(\bar{\k},t)}{\rho_{n+1}(\k_{n+1})\cdot \ldots \cdot \rho_m(\k_m)\mu(t)} 
= \big| \det\big(J_{g_{\bar{\k}}}(t)\big) \big|\, \frac{\prod_{i=1}^n\rho_i(g_{\bar\k,i}(t))}{\mu(t)}.
\end{align*}
		\end{linenomath}
If $A$ is split into subregions, then there is one such expression for each region, with a corresponding density function $\mu(t)$.

\medskip
\paragraph{\bf  Monte-Carlo in practice.}
To approximate the integral $I$ in  \eqref{eq:I}, we sample the $u$ variables from the given distributions,
 and obtain   points  $x^{(i)}=(x_1^{(i)},\dots,x_u^{(i)})$ for $i=1,\dots,N$. 
We then compute $\widehat{I}_{N}$ 
 and the standard error $\hat{e}_N$.  We increase $N$ and sample new points until 
 \begin{equation}\label{eq:stop_criterion}
 \tfrac{\hat{e}_N}{\widehat{I}_N}<10^{-2}. 
 \end{equation} 
 We report  $\widehat{I}_N$ with two digits of significance.
 Some considerations on the minimal sample size are given below.

This method easily allows parallelization.
Specifically, the second expression for the standard error in \eqref{Equation:Standard_error_using_delta} allows for an iterative computation of $\hat{e}_N$ without storing all sampled points, using the cumulative values of $\sum_{i=1}^N Q(x^{(i)})^2$ and $\sum_{i=1}^N Q(x^{(i)})$. 

 \smallskip
As indicated in \cite[\S 2.3]{MonteCarloBook}, the computation of $\hat{e}_N$ using \eqref{Eq:Standard_error_Monte_Carlo} might lead to an imprecise value when $\hat{e}_N$ is very small. 
A way to bypass this problem is to consider $J_1=Q(x^{(1)})$ and $S_1:=0$, 
and iteratively construct the following functions for every new sampled point $x^{(i)}$, $i\geq 2$:
\begin{equation}\label{Equation:Standard_error_using_delta}
\delta_i = Q(x^{(i)}) - J_{i-1},  \qquad 
J_i= J_{i-1} + \tfrac{1}{i} \delta_i,\qquad 
S_i=S_{i-1}  + \tfrac{i-1}{i} \delta_i^2.
\end{equation}
An easy computation shows that $\widehat{I}_N=J_N$ and 
$\hat{e}_N= \sqrt{\frac{S_N}{N(N-1)}}$ (see \cite[\S 2.3]{MonteCarloBook}).

Note that division by $N(N-1)$ for $N$ large may also cause numerical errors. Hence $\hat{e}_N$ is computed by first dividing $S_N$ by $N$, and then by $N-1$. 

\smallskip
When the sample size is too small, then $\widehat{I}_N$ might be an imprecise approximation of the integral $I$, even if the standard error is small. 
This happens when the integrand $f(x)$ in \eqref{eq:I} is nearly zero on $M\setminus M'$,
and the density $p(x)$ of the chosen probability distribution is small on $M'$. 
Then  \eqref{Eq:Integral_Monte_Carlo} and \eqref{Eq:Standard_error_Monte_Carlo} are both close to $0$
if the sample size is too small to cover $M'$ properly (as $f(x)$ will be close to zero for most sampled points) \cite[Ch 9]{MonteCarloBook}. 
 
In practice, it may be difficult to choose the ``best'' probability function.
In this work, we adopt the following thumb rule for the minimum sample size.  
We compute $\widehat{I}_N$ and $\hat{e}_N$ for $N=10^d$, starting with $d=1$. We increase $d$ until $\widehat{I}_N$ belongs to a reasonable interval. For example, if $I$ is the Kac-Rice integral 
of a polynomial system  that \emph{we know}  has between $1$ and $3$ solutions in $A$, then we should have $\widehat{I}_N\in [1,3]$. 
After this initial check on minimum sample size, we  consider the termination condition \eqref{eq:stop_criterion}.

\medskip
\paragraph{\bf Antithetic Monte Carlo. }
When the probability density function is symmetric, then one might 
use antithetic Monte Carlo \cite[\S 8.2]{MonteCarloBook}. 
Specifically, for our setting, consider $u$ independent uniform or truncated normal distributions 
on intervals $[a_i,b_i]$, such that their mean are the centers of  the respective intervals. 
Let $c=(\tfrac{a_1+b_1}{2},\cdots,\tfrac{a_{u}+b_{u}}{2})$ be the center of the product of intervals and $\rho$ the probability density function. Then $\rho(x)=\rho(2c-x)$.   

Antithetic Monte Carlo  consists in sampling $\frac{N}{2}$ points (for $N$ even) and evaluating the   function of interest in each sampled point and its reflection. 
 If computing the reflection of a point is faster than sampling a point, then antithetic Monte Carlo is faster than simple Monte Carlo. Besides, the standard error of antithetic Monte Carlo is between $0$ and $\sqrt{2}\hat{e}$, if $\hat{e}$ is the standard error of simple Monte Carlo \cite[\S 8.2]{MonteCarloBook}.

\medskip
\paragraph{\bf Implementation.}
In our computations, we considered simple and antithetic Monte Carlo implemented manually on different platforms\footnote{Versions: \texttt{Maple 2020}, \texttt{Python 3.7.4}, \texttt{C++11}, \texttt{Numba 0.48.0} and \texttt{Julia 1.4.2}.}: \texttt{Maple}, \texttt{Python} (with and without the package \texttt{Numba}), \texttt{C++} and \texttt{Julia}. Additionally, we considered the \texttt{CUBA} package \cite{CUBA-origin}, as already implemented in all these platforms\footnote{List of platforms providing the CUBA package: \url{http://www.feynarts.de/cuba/}.}, which has four advanced numerical integration techniques: Vegas (Monte Carlo integration with importance sampling), Suave (Monte Carlo integration with globally adaptive subdivision and importance sampling), Divonne (Monte Carlo integration with stratified sampling and numerical optimization), and Cuhre (deterministic integration with globally adaptive subdivision).

Despite the advanced techniques implemented in \texttt{CUBA}, the computation of the Kac-Rice integral in \texttt{CUBA}  often gives inaccurate answers in examples, while a manual implementation works well. The problem with CUBA persisted even after changing accuracy levels and options values for all four available integration methods. 
We attribute the problem to the fact that we cannot choose the distribution for sampling in \texttt{CUBA}. This observation additionally supports the appropriateness of our choice of distribution.

The integrals reported in Section~\ref{sec:examples} have been computed using our own implementation\footnote{Computations performed per default on Windows 10, Intel(R) Core(TM) i7-2670QM CPU @ 2.20GHz 2.20 GHz,  x64-based processor, 6.00GB (RAM)} with \texttt{Julia}. In Subsection~\ref{Example:Joshi_line_2v5p_network_antithetic_uniform}, we compare the speed of computation of the Kac-Rice integral for a specific example using \texttt{Maple}, \texttt{Python} (with and without \texttt{Numba}), \texttt{C++} and \texttt{Julia}.  The analysis showed that \texttt{Numba} was the fastest option, competing closely with  \texttt{Julia}. However, manual parallelization using \texttt{Julia} and the package \texttt{Distributed} is easier than with \texttt{Numba} using the module \texttt{multiprocessing}. Therefore we favoured \texttt{Julia} over \texttt{Numba}. For testing parallelization, we have used a server consisting of 64 cpus, AMD Opteron(tm) Processor 6380.

The code for the computations in Section~\ref{sec:examples} can be found in a GitHub repository archived by Zenodo \cite{MCKR_project_repository}. A separate folder contains the relevant files for each subsection.
 A \texttt{Julia} program named MCKR and available at  \cite{MCKR_program} can be used to apply the methods in next section to any example that satisfies the assumptions of Theorem~\ref{thm:Kac-Rice}. The user needs to provide (in \texttt{Julia})
the functions $g_i$, the determinant of the Jacobian, the choice of random distribution, and the desired task. See the manual available at \cite{MCKR_program} for information.

\section{Parameter regions using Kac-Rice formulas} \label{sec:examples}
As stated in the introduction and motivated in Subsection~\ref{sec:Reaction_networks_and_multistationarity}, our main goal is to understand the parameter region in relation to the number of  positive solutions of a polynomial system. For
a parametrized polynomial system $f_\k(t)=0$, we focus on determining the expected number of positive solutions when the parameters $\k$ belong to a bounded box $B$. To this end, we compute the Kac-Rice formula after  endowing all parameters with  uniform distributions. That is, let $B=B_1\times \cdots \times B_m$ with $B_i$ bounded intervals, consider
$\k_i \sim U(B_i)$ for $i=1,\dots,m$, and let 
\[ \hat{r}(B)=\mathbb{E}(f_{\k}^{-1}(0) \cap \R^n_{>0}).\]
Then $\hat{r}(B)$ is the average number of positive solutions of the system $f_\k(t)=0$ for $\k\in B$.

Let $M_{\rm max}$ and $M_{\rm min}$ be the maximal and minimal number of positive solutions the system $f_\k(t)=0$ generically admits (that is, in some open set of $\R^m$). If $\hat{r}(B)= M_{\rm min}$ resp. $ M_{\rm max}$, then for almost all parameter values in $B$, the system has 
$M_{\rm min}$, resp.  $M_{\rm max}$ solutions. 
If $M_{\rm min}<\hat{r}(B)< M_{\rm max}$, then all we can assert is that $B$ contains 
parameter values where the system has more than $M_{\rm min}$ solutions. In general, if $\hat{r}(B)>M$ for some $M$, then $B$ contains parameter points $\k$ where $f_\k(t)=0$ has more than $M$ positive solutions.

We aim at  dividing the parameter region into areas where 
\begin{equation}\label{eq:classifyI}
{\rm (i)} \ \hat{r}(B)=M_{\rm max},\quad {\rm (ii)} \ \hat{r}(B)=M_{\rm min}, \quad \textrm{and}\quad {\rm (iii)} \  M_{\rm min}<\hat{r}(B)< M_{\rm max}, 
\end{equation}
or in the setting of reaction networks, into areas where 
\begin{equation}\label{eq:classifyII} {\rm (i)} \ \hat{r}(B)=M_{\rm max},\quad {\rm (ii)} \ \hat{r}(B)\leq 1, \quad \textrm{and}\quad {\rm (iii)} \  1<\hat{r}(B)<M_{\rm max}. 
\end{equation}
In this scenario, cases (i) and (iii) include the region of multistationarity if $M_{\rm max}>1$. 

With this in mind, we use $\hat{r}(B)$
 (if it is well defined and can be computed) to address the following two problems.

\medskip
\paragraph{\bf Problem I: Coarse description of parameter regions of multistationarity.} 
 Let $\delta_1,\dots,\delta_m>0$ be the desired precision for each parameter, that is, the minimal lengths of the intervals $B_1,\dots,B_m$ to consider. Consider a grid partition of some box $B=B_1\times \dots \times B_m$ in small sub-boxes $C_1,\dots,C_\ell$ of side length at most $\delta_i$ for the $i$-th variable. We approximate the classification of the parameter points according to the number of solutions to $f_\k(t)=0$ by computing $\hat{r}(C_i)$ for $i=1,\dots,\ell$ and classifying it into cases (i)-(iii) as in \eqref{eq:classifyI} or \eqref{eq:classifyII}. 
 In the setting of reaction networks, using \eqref{eq:classifyII} we obtain a coarse approximation of the real parameter region of multistationarity, as well as the region where multistationarity does not occur.

In order to optimize the speed of computation, we use a \textbf{bisect strategy}.   
 If $\hat{r}(B)\neq  M_{\rm min}, M_{\rm max}$, then we bisect $B$ along one direction, and obtain two sub-boxes $C_1,C_2$. We compute $\hat{r}(C_1)$ and $\hat{r}(C_2)$. 
If $C_1$ belongs to cases (i) or (ii) of \eqref{eq:classifyI} or \eqref{eq:classifyII}, then we have classified this box and move onto $C_2$. Otherwise, if some side of the box is larger than $\delta_i$, we repeat with $C=C_1$. We perform the same procedure with $C_2$. 

We start by considering the maximal number of steps for each parameter value $\k_i$, as given by the precision $\delta_i$:
\begin{equation}\label{eq:stop_precision}
 L_i:= {\rm ceiling}\Big(\log_{2} \Big( \tfrac{{\rm length}(B_i)}{\delta_i}\Big)\Big), \qquad i=1,\dots,m.
\end{equation}
At the $j$-th step, the direction of bisection is
the axis along the parameter $\k_i$ for which $i=j$ (mod $m$). If the direction of $\k_i$ has already been bisected  $L_i$ times, then this direction is skipped. 

In this way, larger boxes already belonging to (i) or (ii) are not subdivided and hence the computational time is reduced dramatically. 
This approach considers smaller boxes containing the boundary separating regions where the number of solutions to $f_\k(t)=0$ changes. 
 
With a grid description of the parameter region of multistationarity, one can derive a semialgebraic set 
defined by a single polynomial,  which contains the multistationarity region, and with the minimal volume (see \cite{PolynomailSuperlevelSetRepresentation}). Additionally, the boundary of the region of multistationarity can be approximated by the hypersurface given by the polynomial in the superlevel set representation.

\medskip
\paragraph{\bf Problem II: Parameter point with multistationarity. } We aim at finding a parameter point or sub-box for which the system has $M_{\rm max}$ solutions in a  given bounded box $C$ of interest, or conclude that no such parameter choice exists. 
To this end, we apply the bisect strategy, but keeping at each step the sub-box with largest $\hat{r}$, and stopping when $\hat{r}(C) = M_{\rm max}$  (approximately) or the maximal number of divisions has been reached for all parameters.

If the precision is small enough, this strategy  is guaranteed to work if the system only has two possible number of solutions for generic parameter values. If that is not the case, then we might not identify a box with $M_{\rm \max}$ solutions. For example, if the system generically admits one, three or five positive solutions, and at one step the two boxes $C_1$ and $C_2$ to consider are such that $C_2$ belongs to the region with three solutions, while $C_1$ intersects the regions with one and with five but such that $\hat{r}(C_2)>\hat{r}(C_1)$, we will miss the region with five solutions. 
To bypass this problem, we should search the parameter region as in \textbf{Problem I}, and keep both boxes unless $\hat{r}$ equals $M_{\rm min}$.

This  approach can also be used to numerically determine the maximal number of positive solutions the system admits in a box $C$, and to search for parameter points for which the system has a given number of solutions $M$. If $M\neq M_{\rm max}, M_{\rm min}$, then $\hat{r}(C)\simeq M$ does not guarantee that all parameters in the box give rise to $M$ solutions, so one needs to pick a point and verify the number of solutions by solving the system. 
Alternatively, in  \cite[Lemma 5.4]{PolynomailSuperlevelSetRepresentation} it is shown that by considering a distribution on $\k$ different from the uniform,  one can check whether $\hat{r}(C)\simeq M$ implies that all parameters in the box give rise to $M$ solutions.

Finally, if of interest is only to determine the existence of parameter values for which the system has more than one solution, then it is enough to find a box $C$ with $\hat{r}(C)>1$. Then for any parameter $\k$ in $C$, 
the box $\prod_{i=1}^m [\k_i-\delta_i,\k_i+\delta_i]$, contains a point where the system has more than one solution.

 \smallskip
 Theoretically,  these two problems  can be addressed using CAD \cite{basu2007algorithms,FabriceUsingCAD,SolvingParametricPolynomialSystemsFabrice}. However, this method is impractical as it is double exponential in the total number of variables and parameters, and depends also on the number and degree of the polynomials   \cite{NumberOfCellsCAD}.
There are theoretically faster algorithms based on the \textit{critical points method}, which returns a finite set of points including at least one point from each connected component of a semi-algebraic set \cite{basu2007algorithms}. This method, of single exponential complexity in the number of variables, can be used to address \textbf{Problem II}, by considering the semi-algebraic set given by the defining inequalities of the box $B$, together with $p\neq 0$, where $p$ is the polynomial defining the discriminant variety of the parametric system~\eqref{eq:steadystates}. The number of solutions of system~\eqref{eq:steadystates} is invariant in each connected component of this semi-algebraic set. 
  When combined with roadmap algorithms, that decide whether two points belong to the same connected component, the number of connected components and even semi-algebraic descriptions of the components can be found \cite[Chapters 15\&16]{basu2007algorithms}. Using this approach, \textbf{Problem I} can also be addressed. 
An algorithm to study \textbf{Problems I} and \textbf{II} using these ideas is singly exponential in the number of parameters of the system~\eqref{eq:steadystates}, and doubly or singly exponential in the number of variables of the system~\eqref{eq:steadystates}. The later is a consequence of the computation of the discriminant variety of a parametric system. To this end, there are several possible approaches. The most common approach is to use elimination theory via Gr\"obner basis computation. Gr\"obner basis computation is known to be doubly exponential in the worst case \cite{Mayr-Meyer-GrobnerBasisComplexity,Mayr-Ritscher-GrobnerBasisComplexity}. An alternative approach is to use the projection operator of CAD algorithms with respect to the variables of the system only  \cite{Amir-Matthew-ImprovingAlgebraicTools}. The projection step of CAD is still doubly exponential on the number of variables \cite{Bradford-Davenport-England-McCallum-Wilson-DoublyExponentialityProjection}. However, using Equational Constraints, the complexity can be reduced to singly exponential on the number of variables  \cite{Amir-Matthew-ImprovingAlgebraicTools,England-Bradford-Davenport-Equational-Constraints}.

In what follows we provide several examples (mainly arising from reaction networks) to illustrate how to address the two problems described above by computing $\hat{r}(C)$ using the Kac-Rice formula. For small examples, we compare our results with the exact answer.
As an effective implementation for the critical points method approach mentioned above is yet to be developed, we compare our results with the output of CAD using the package 
\texttt{RootFinding[Parametric]}  of \texttt{Maple 2020} \cite{MaplePackageCAD}. 
 
  We start in Subsection~\ref{Example:HK_network} with an illustrative reaction network with eight parameters where the number of positive steady states is generically one or three, and the system \eqref{eq:steadystates} can be reduced to one polynomial equation. For illustration purposes, we start by fixing the value of six parameters, finding  the parameter regions of interest, 
 and comparing them visually to the output of CAD.  Afterwards,  we show that {\bf Problems I and II} can also be solved with eight free parameters.  
 
  We proceed with another reaction network in Subsection~\ref{Example:Joshi_line_2v5p_network_antithetic_uniform} with five free parameters. 
 We find a parameter point with multistationarity and  compare the performance of simple and antithetic Monte Carlo in  different platforms. 
 
 We next study a polynomial in one variable and two parameters that admits five positive roots (Subsection~\ref{Example:Degree5}).  We study the partition of the parameter space according to the number of positive roots of the polynomial as given by the Kac-Rice formula and compare the result with CAD. 
 
 Finally, we study two relevant reaction networks in Subsection~\ref{sec:twosite} and \ref{sec:extendedHK}, with a higher number of parameters and variables.

\subsection{Illustrative example: two component system}\label{Example:HK_network}
The following reactions define a reaction network representing a simplified model of a  two-component system with hybrid histidine kinase as considered in \cite{FeliuHK}:
\begin{equation}\label{eq:hybridHK}
\begin{aligned}
	X_1  \ce{->[k_1]} X_2 &  \ce{->[k_2]} X_3 \ce{->[k_3]} X_4 \quad & 
	X_3+X_5 &  \ce{->[k_4]} X_1+X_6\\
	X_6&  \ce{->[k_6]} X_5 & X_4+X_5 & \ce{->[k_5]} X_2+X_6.
\end{aligned}
\end{equation}
The system of parametrized polynomial equations \eqref{eq:steadystates} is
\begin{align*}
	k_4x_3x_5-k_1x_1 & =0, &
	k_5x_4x_5+k_1x_1-k_2x_2 & =0,\\
	-k_4x_3x_5+k_2x_2-k_3x_3& =0,&
	-k_4x_3x_5-k_5x_4x_5+k_6x_6& =0,\\
	x_1+x_2+x_3+x_4-T_1& =0,&
	x_5+x_6-T_2& =0.
\end{align*}
It is shown in \cite{FeliuHK} 
that the positive solutions to this system are in one-to-one correspondence with the positive solutions to the following univariate polynomial of degree three in $t=x_5$:
\begin{equation}\label{Equation:Univariate_polynomial_HK_8_parameters_form_t} 
	\begin{split}
	f_\k(t)=(k_1+k_2)k_4k_5k_6t^3+(T_1k_1k_2k_4k_5-T_2(k_1+k_2)k_4k_5k_6+	k_1(k_2+k_3) k_5k_6)t^2 \\ 
	+(T_1k_1k_2k_3k_5-T_2k_1k_5k_6(k_2+k_3)  	+k_1k_2k_3k_6)t-T_2k_1k_2k_3k_6. 
	\end{split}
\end{equation}
Thus, in this example, the goal is to study the number of positive roots of a degree three polynomial, as function of the eight parameters $k_1,\dots,k_6>0$ and $T_1,T_2$.  
As shown in \cite{FeliuHK}, there exist parameter values for which \eqref{Equation:Univariate_polynomial_HK_8_parameters_form_t}  has three positive roots, and it always has at least one. 		Observe that we necessarily have $T_1,T_2>0$ for positive solutions to exist.
CAD is computationally prohibitive with $8$ parameters $k_1,\dots,k_6,T_1,T_2$ on a standard computer. As we will see below, the Kac-Rice formula combined with Monte Carlo integration can cope with this situation. 

\begin{figure}[t]
		\centering
		\begin{minipage}[h]{0.3\textwidth}
				\begin{center}
				\includegraphics[width=4cm]{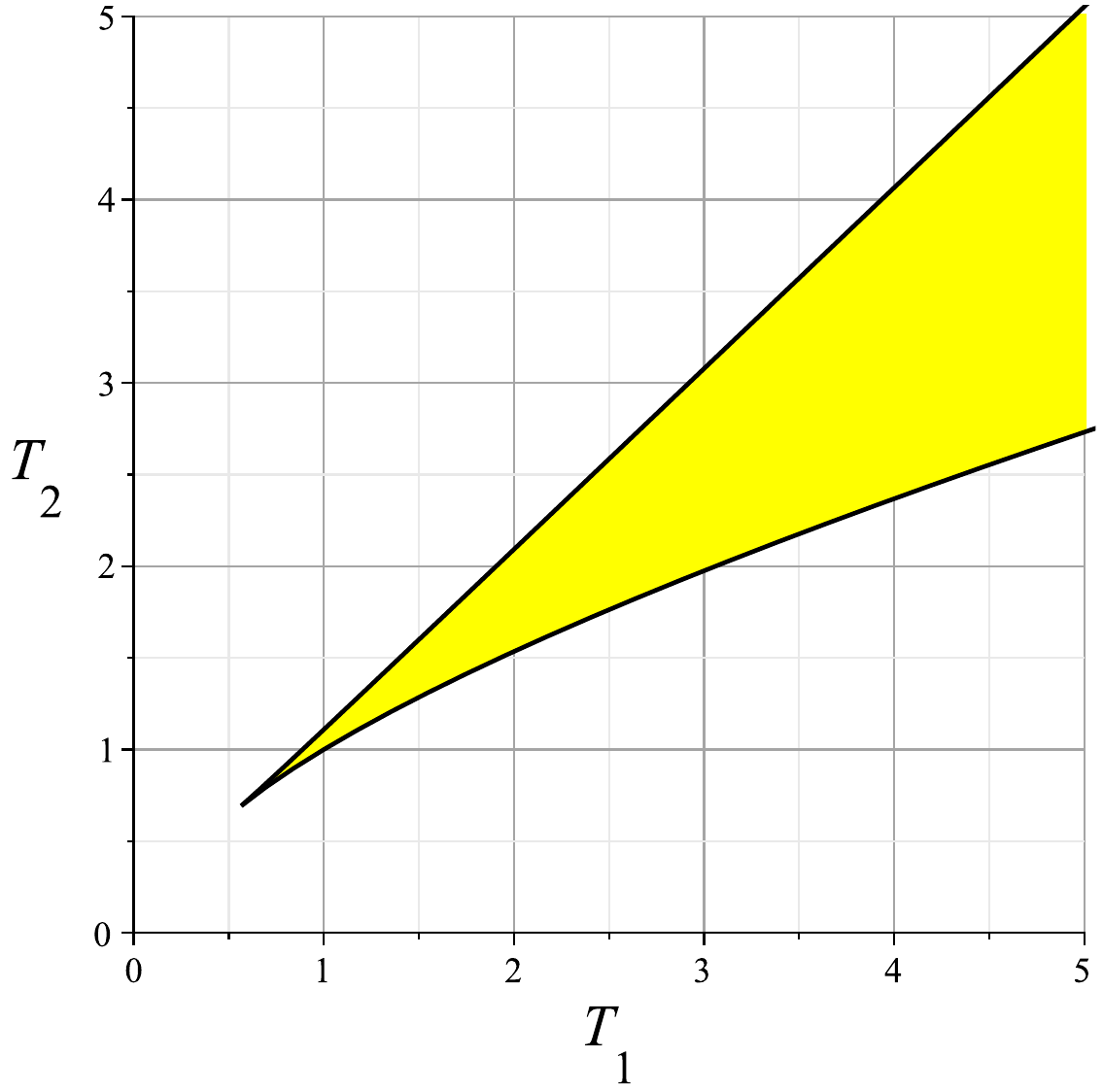}
				
				(a)
				\end{center}
		\end{minipage}
		\begin{minipage}[h]{0.3\textwidth}
				\begin{center}
				\includegraphics[height=4cm]{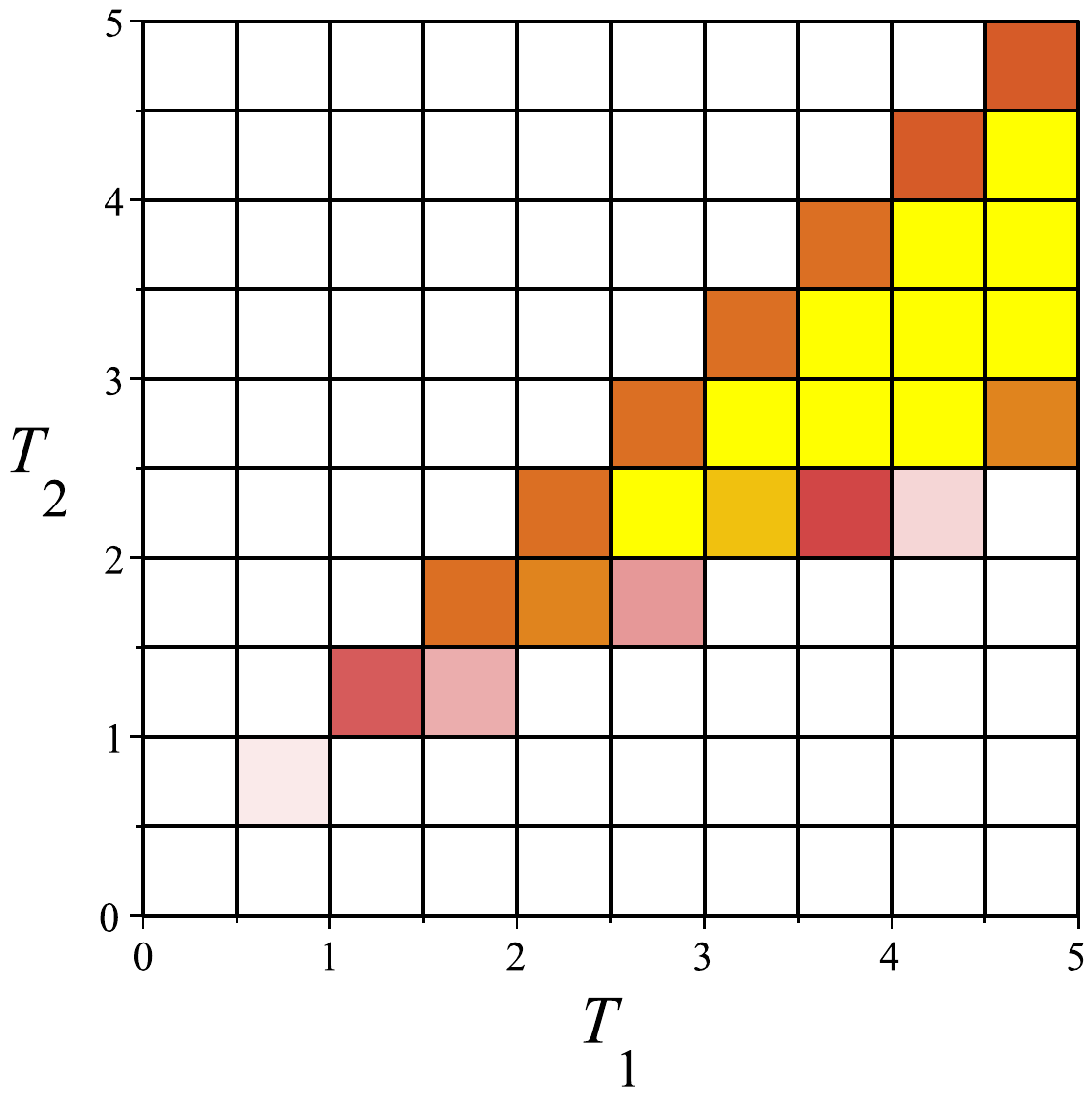}
				
				(b)
				\end{center}
		\end{minipage}
		\begin{minipage}[h]{0.06\textwidth}
			\begin{center}
				\vspace{-1cm}\includegraphics[height=3.5cm]{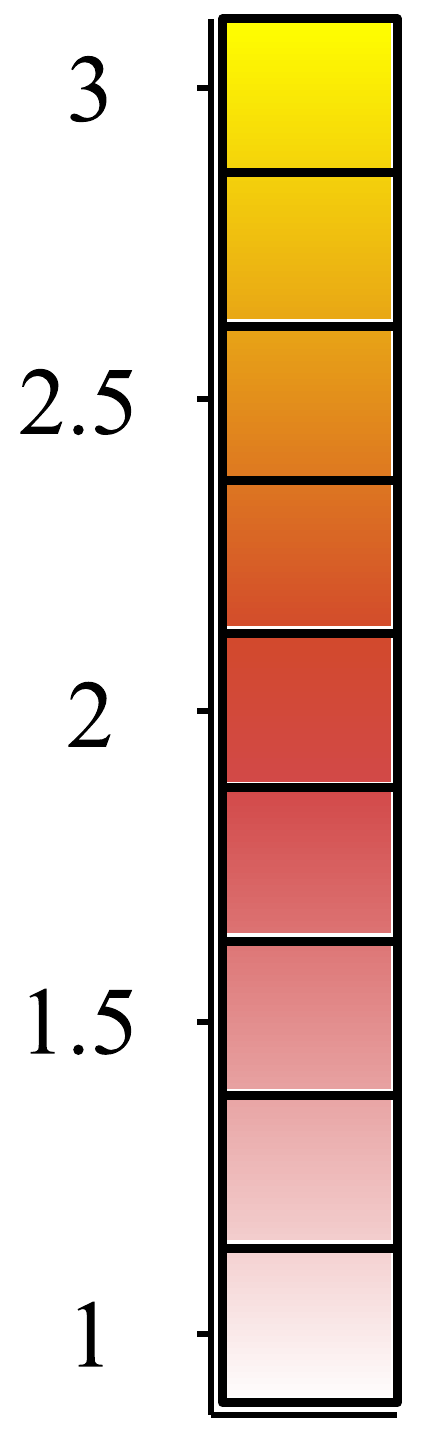}

			\end{center}
		\end{minipage}
		\begin{minipage}[h]{0.3\textwidth}
			\begin{center}
			\includegraphics[height=4cm]{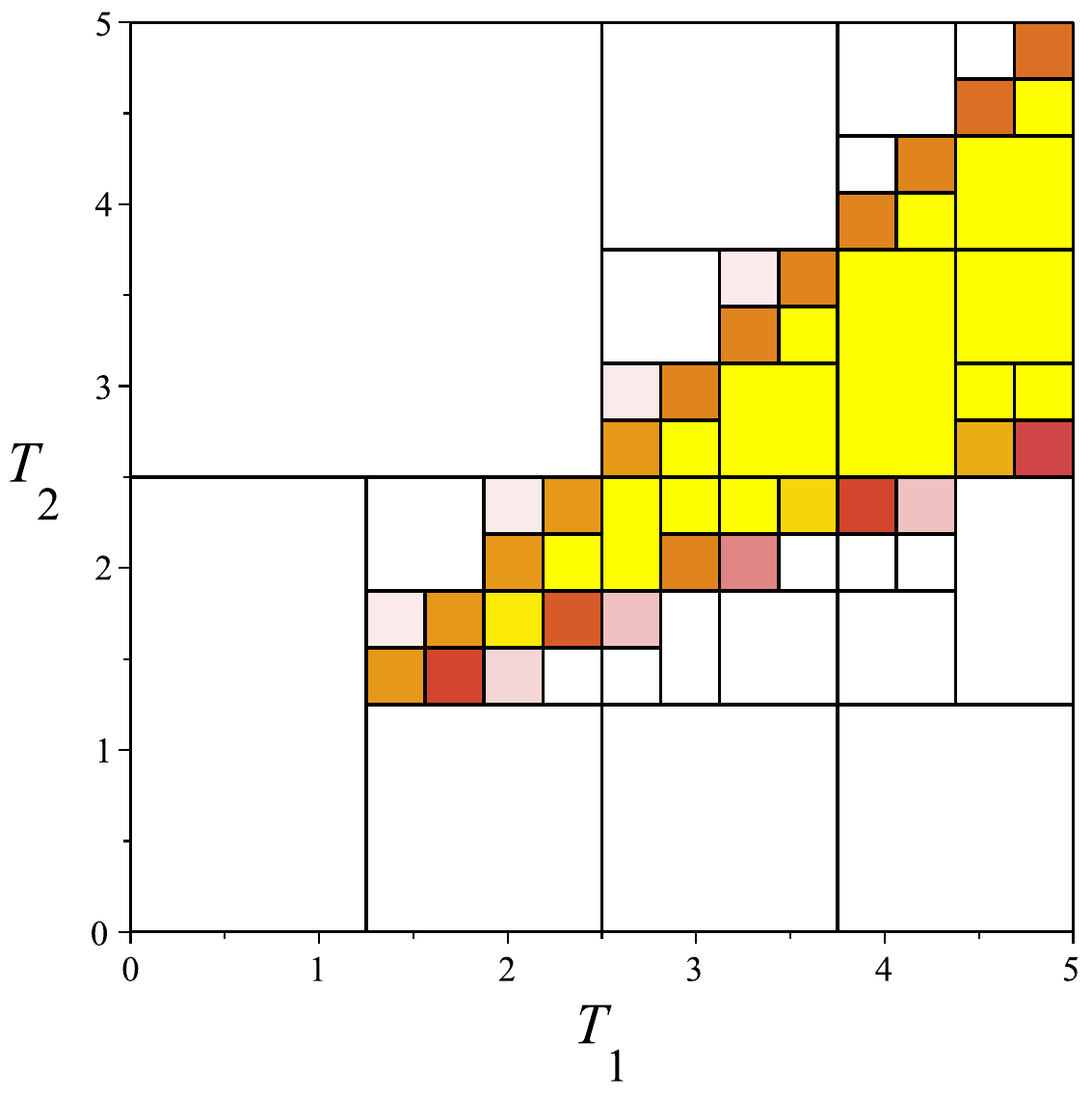}
			
			(c)
			\end{center}
		\end{minipage}
		\caption{{\small Parameter regions in $T_1,T_2$  according to the number of positive steady states for  network  \eqref{eq:hybridHK} with $k$ as  in \eqref{Equation:Reaction_rate_constant_choice_HK_2_parameters}.  (a) Obtained using CAD. The network  has three positive steady states in the yellow region and one in the white region (and two steady states on the boundary between the two regions, one with multiplicity two).
		(b-c) Approximation of the parameter region using the Kac-Rice formula and numerical integration on sub-boxes.  Yellow corresponds to three positive steady states and white to one, see bar diagram. For (b) the sub-boxes are constructed by grid partitioning, while for (c) the sub-boxes are constructed by the bisect strategy with termination condition of $4$ divisions per parameter.}} 
		\label{Figure:HK_CAD}	\label{Figure:HK_Uniform}
	\end{figure}

\medskip
\paragraph{\bf Identifying the region of multistationarity. } 
For \emph{illustrative purposes}, we first fix the  reaction rate constants $k_i$ and understand the region defined by the parameters $(T_1,T_2)$ according to the number of positive roots of the polynomial. 
	In \cite{FeliuPlos} 
	it is shown that there exists a choice of $(T_1,T_2)\in\R^2_{>0}$ for which the network is multistationary if and only if $k_1<k_3$.
	So we fix the following reaction rate constants (from \cite[Fig.~2C]{FeliuHK}):
	\begin{linenomath}
	\begin{equation}\label{Equation:Reaction_rate_constant_choice_HK_2_parameters}
	(k_1,\dots,k_6)=(0.7329,100,73.29,50,100,5).
	\end{equation}
	\end{linenomath}
	Evaluating the univariate polynomial \eqref{Equation:Univariate_polynomial_HK_8_parameters_form_t} at \eqref{Equation:Reaction_rate_constant_choice_HK_2_parameters} gives a polynomial $f_{T_1,T_2}(t)$ of degree 3 in $t$, whose coefficients depend on the two parameters $T_1$ and $T_2$:
\begin{linenomath}	\begin{align}\label{eq:f_evaluated}
	\begin{split} 
	f_{T_1,T_2}(t)= & (2518322.5)t^3+\big((366450)T_1-(2518322.5)T_2+63502.1205\big)t^2 \\
	& +\big((537142.41)T_1-(63502.1205)T_2+26857.1205\big)t-(26857.1205)T_2.
	\end{split}
	\end{align}\end{linenomath}
	The analysis of this polynomial is addressable using CAD, which provides an explicit description of the region where \eqref{eq:f_evaluated} has three positive solutions. For $T_1,T_2\in (0,5)$, the region is depicted in Figure~\ref{Figure:HK_CAD}(a).
	
We consider now the same problem using the Kac-Rice formula and Monte Carlo integration. In the notation of Theorem~\ref{thm:Kac-Rice}, by letting $\bar{\k}=T_2$ and $n=1$, we have
\begin{linenomath}\begin{align*}
h(t) &= 366450t^2+537142.41t, \\
q(T_2,t)& =-(2518322.5t^2+63502.1205t+26857.1205) (T_2-t).
\end{align*}\end{linenomath}
 Then, for any bounded box $B= [a_7,b_7]\times [a_8,b_8]$,  we have
\begin{linenomath}\begin{align}
\hat{r}(B)=\int_0^{+\infty}   \int_{a_8}^{b_8} \frac{ |J_{T_2}(t)|}{(b_8-a_8)(b_7-a_7)} \chi_{[a_7,b_7]}\big(\tfrac{-q(T_2,t)}{h(t)}\big)dT_2\, dt,\label{Equation:Kac-Rice_integral_HK_1v2p}
\end{align}\end{linenomath}
where $J_{T_2}(t) = \tfrac{\partial}{\partial t}\Big(\tfrac{-q(T_2,t)}{h(t)}\Big)$. 
As $t=x_5$ and $0<x_5,x_6$,  $x_5+x_6=T_2$, any positive root of  \eqref{eq:f_evaluated} for parameter values in $B$ satisfies $t<b_8$. 
Hence we choose $\mu(t)$ (see Subsection~\ref{sec:Numerical_integration}) to be the density of the uniform distribution on $(0,b_8)$. 
The Kac-Rice integral is then approximated by the following sum for randomly sampled points $t^{(i)},T_2^{(i)}$ for $i=1,\dots,N$ and $N$ large:
\begin{equation}\label{eq:example_HK_1v2p_MonteCarlo}
\tfrac{b_8}{(b_7-a_7)}\sum_{i=1}^N\big| J_{T_2^{(i)}}(t^{(i)}) \big| \chi_{[a_7,b_7]}\Big(\tfrac{-q(T_2^{(i)},t^{(i)})}{h(t^{(i)})}\Big).
\end{equation}

We consider the box $B=[0,5]\times [0,5]$, subdivide it into 100 sub-boxes (of side length $0.5$), and for each sub-box compute $\hat{r}(B)$ using \eqref{eq:example_HK_1v2p_MonteCarlo}. It took 46 seconds and 100 integrals were computed.
We depict the output in Figure~\ref{Figure:HK_Uniform}(b), where we color each sub-box with a graduation of yellow, orange and white:  yellow means the expected number is three, and white means it is one.

Clearly, Figure~\ref{Figure:HK_Uniform}(b) approximates Figure~\ref{Figure:HK_CAD}(a), which displays the exact region.  In Figure \ref{Figure:HK_Uniform}(b) the sub-boxes that cross the thick line separating the yellow and white regions in Figure~\ref{Figure:HK_CAD}(a) have an orange-like color, because the sub-box contains parameters with both one and three positive steady states. By making the size of the sub-boxes smaller, we would get more accurate approximations of Figure~\ref{Figure:HK_CAD}(a). 

Figure~\ref{Figure:HK_Uniform}(c) has been found using the bisect strategy.  
For the minimal box length to be at most $0.5$, \eqref{eq:stop_precision} gives that $4$ bisections are (at most) required for each parameter. The process took 52 seconds, computed 111 integrals and returned 56 sub-boxes.

\smallskip
 
This example  illustrates how the Kac-Rice formula can be used to approximate the parameter region. 
The advantage is that the numerical integrals we need to compute require, in principle, less computer power than performing CAD.

\medskip
\paragraph{\bf Finding a multistationary point in a box. }
We consider now the problem of finding a parameter value where \eqref{Equation:Univariate_polynomial_HK_8_parameters_form_t} has three positive roots. 
We follow the approach outlined for {\bf Problem II} at the beginning of this section.

	\begin{figure}[t!]
	\begin{minipage}[h]{0.55\textwidth}
	
		\centering
		\begin{tabular}{|c|c|c|c|}
			\hline
			Step & Sub-box $B$ & $\hat{r}(B)$ & \begin{minipage}[h]{0.2\textwidth} {\small Chosen\\[-1pt] sub-box \\[-10pt] }\end{minipage} \\ \hline
			0 & $[1,3]\times [2,4]$ & $\simeq 1.29$ & \checkmark \\ \hline
			1 & $[1,2]\times [2,4]$ & $\simeq 1.00$ &            \\ \cline{2-4} 
			& $[2,3]\times [2,4]$ & $\simeq 1.58$ & \checkmark   \\ \hline
			2 & $[2,3]\times [2,3]$ & $\simeq 2.16$ & \checkmark \\ \cline{2-4} 
			& $[2,3]\times [3,4]$ & $\simeq 1.00$ &              \\ \hline
			3 & $[2,2.5]\times [2,3]$ & $\simeq 1.68$ &          \\ \cline{2-4} 
			& $[2.5,3]\times [2,3]$ & $\simeq 2.65$ & \checkmark \\ \hline
			4 & $[2.5,3]\times [2,2.5]$ \cellcolor{yellow} & $\simeq 3.00$\cellcolor{yellow}  & \checkmark \cellcolor{yellow} \\ \cline{2-4} 
			& $[2.5,3]\times [2.5,3]$ & $\simeq 2.30$ &          \\ \hline
		\end{tabular}
		
		\medskip
		(a)
	\end{minipage}
	\quad
 	\begin{minipage}[h]{0.35\textwidth}
	\centering
		\includegraphics[height=5cm]{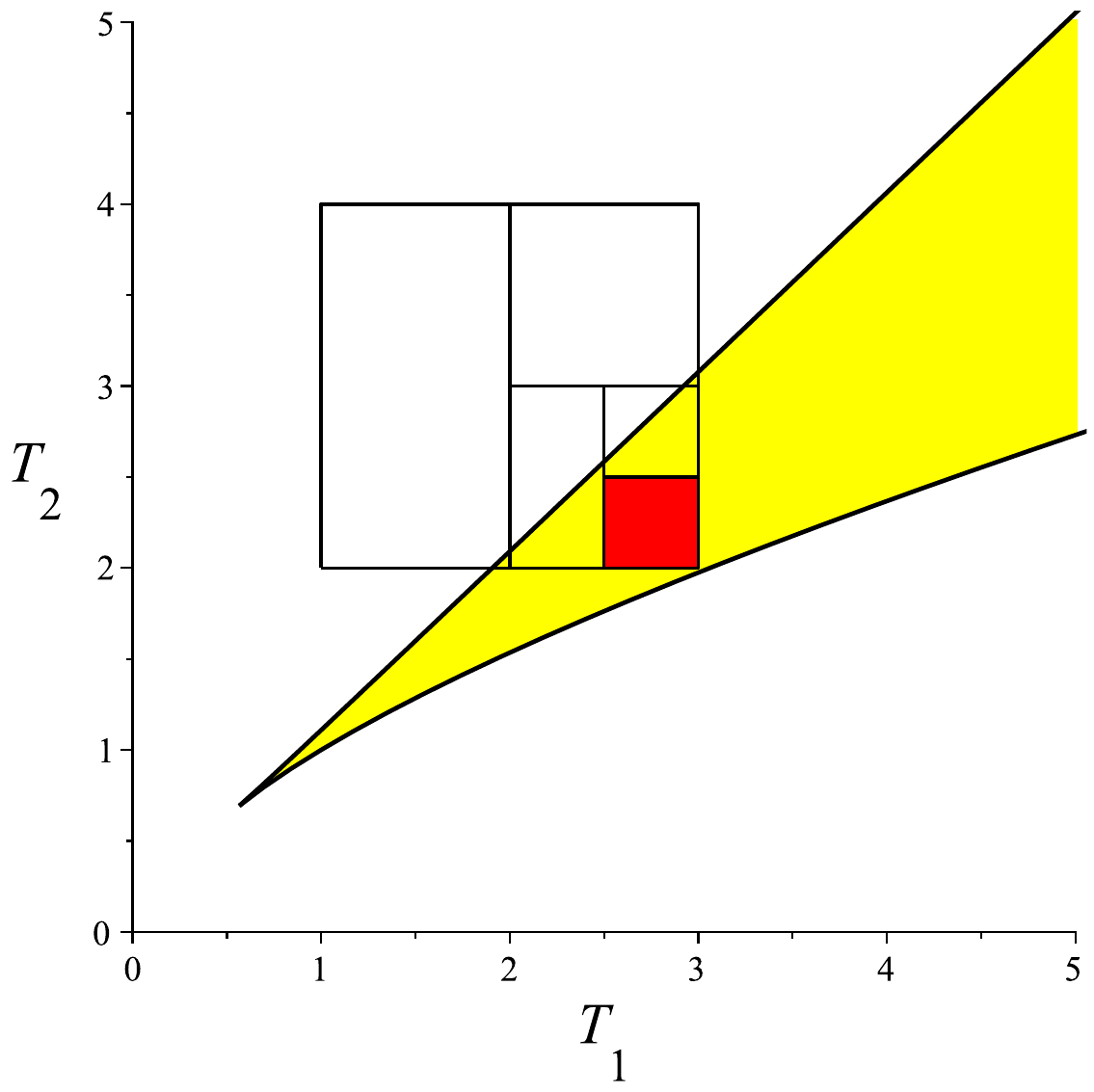}
		
		(b)
	\end{minipage}
		\caption{{\small {\bf Problem II} for network \eqref{eq:hybridHK} with $k$ as in \eqref{Equation:Reaction_rate_constant_choice_HK_2_parameters}. At each step, $\hat{r}$ is computed for the considered sub-boxes,  and the one with the largest $\hat{r}$ is bisected, until $\hat{r}$ is  $3$ with two decimal digits of precision. At step $0$, the test verifies that $\hat{r}>1$, otherwise the box does not intersect the region of multistationarity.  Here the $T_1$-axis is chosen in odd steps, and the $T_2$-axis in even steps.  (a) Table description of the considered sub-boxes and their $\hat{r}$. (b) Visual depiction of the sub-boxes in (a). The final sub-box is colored in red and is entirely inside the multistationary region (colored in yellow). The sub-boxes with $\hat{r}=1$ are outside the yellow region, and the sub-boxes with $1<\hat{r}<3$, have intersection with both the white and yellow regions.
}}
		\label{Figure:Search_box}
	\end{figure}

	Again for  illustrative purposes, fix the reaction rate constants in \eqref{Equation:Reaction_rate_constant_choice_HK_2_parameters} and consider the box $B=[1,3]\times [2,4]$
	in  the parameter space for  $T_1$ and $T_2$.
	Computing the Kac-Rice integral \eqref{Equation:Kac-Rice_integral_HK_1v2p} we find $\hat{r}(B)\simeq 1.29$, and hence there are parameter values in $B$ yielding more than one positive steady state.
We proceed to iteratively bisect $B$ and compute the Kac-Rice integral for the two resulting sub-boxes, until we obtain a sub-box $B'$ with $\hat{r}(B')\simeq 3$. 
	  
	  Figure~\ref{Figure:Search_box}(a)  shows the result of an implementation of this process, and Figure~\ref{Figure:Search_box}(b) depicts the sub-boxes considered in the process and highlights the found sub-box. Figure~\ref{Figure:Search_box}(b)  shows the real region with three positive steady states from Figure~\ref{Figure:HK_CAD}(a) in the background, such that one can visually interpret the expected number of solutions given in  Figure~\ref{Figure:Search_box}(a), and verify that the final sub-box is entirely inside of the multistationarity region.

\medskip
\paragraph{\bf With $8$ parameters. } In the previous analysis of network~\eqref{eq:hybridHK}, we kept 
only 2 parameters free to be able to visually illustrate our approach, as well as to compare with CAD. We  show here that we can find a parameter point/box where the network has three positive steady states, also when all 8 parameters are free.

We consider the following  box for the parameter vector $(k_1,k_2,k_3,k_4,k_5,k_6,T_1,T_2)$:
\[ B = (0,1)\times (0,200) \times (0,100) \times (0,100) \times (0,200) \times (0,10) \times (0,5) \times (0,5). \]
Computing the Kac-Rice integral we find $\hat{r}(B)\simeq 1.2$ with two digits of significant. Therefore it has intersection with the multistationarity region. We apply the algorithm for \textbf{Problem II}, which, after 22 iterations in less than 58 seconds,
returns the following sub-box:
\begin{linenomath}
\begin{multline*}
C = (0.125,0.25)\times (125,150) \times (75,87.5) \times (12.5,25) \times (175,200)\\ \times (2.5,3.75) \times (3.75,5) \times (3.75,5). 
\end{multline*}\end{linenomath}
For almost all parameter values in this box, the network has three positive steady states, because $\hat{r}(C)=2.97$ with  standard error $\hat{e}=0.009$.

\medskip
We address also {\bf Problem I} with the bisect strategy to obtain a coarse description of the parameter region of multistationarity inside the following box $B$,
\[(0.125,0.375)\times (100,125)\times (75,100)\times (12.5,37.5)\times (150,200)\times (1.25,3.75)\times (0,5)\times (0,5). \]
As input precision, we considered $\delta=(\delta_1,\dots,\delta_8)$ (the upper bound of the minimal length of the intervals  for $\k_1,\dots,\k_8$) as follows: 
\[\delta = (0.125,\,12.5,\,12.5,\,12.5,\,25,\,1.25,\,1.25,\,1.25).\]
The algorithm terminates after computing $1787$ integrals in $2164$ seconds. The obtained partition consists of $894$ sub-boxes, $204$ of which are out of the multistationary region, $334$ are inside the multistationary region, and the remaining $356$ have intersection with both the region of multistationarity and monostationarity.

\subsection{Method and platform comparison: An example with $5$ parameters. }\label{Example:Joshi_line_2v5p_network_antithetic_uniform}
Consider the following reaction network 
	\begin{equation}\label{Equation:Joshi_2v5p_network}
	\xymatrix @C=1.5pc @R=1pc {2X_1+X_2\ar[r]^{\quad k_1} & 3X_1\ar[r]^{k_2\quad} & X_1+2X_2\ar[r]^{\quad k_3} & 3X_2\ar@/^1.0pc/[lll]^{\quad k_4}}.
	\end{equation}
System~\eqref{eq:steadystates} becomes a parametrized polynomial system in five parameters and two variables: 
	\begin{align}
	\k_1t_1^2t_2-2\k_2t_1^3-\k_3t_1t_2^2+2\k_4t_2^3& =0, & 
	t_1+t_2-\k_5& =0,
	\end{align}
	where $\k_i=k_i>0$, $i=1,\dots,4$ and $\k_5>0$. 
We find $\hat{r}(B)$ for $B=\prod_{i=1}^5(a_i,b_i)$ a box in the parameter space with  $a_i\geq 0$. 
In order to apply Theorem \ref{thm:Kac-Rice}, we choose $\k_1$ and $\k_5$ as the linear parameters, which gives $\bar{\k}=(\k_2,\k_3,\k_4)$,
	\begin{align*}
g_{\bar{\k},1}(t) &=\tfrac{1}{t_1^2t_2}(2\k_2t_1^3+\k_3t_1t_2^2-2\k_4t_2^3),&
g_{\bar{\k},2}(t) &=t_1+t_2,
	\end{align*}
and 
\[ \det( J_{g_{\bar\k}}(t))= \tfrac{1}{t_1^3t_2^2} (t_1+t_2)(2\k_2t_1^3-\k_3t_1t_2^2+4\k_4t_2^4).\]
The numerator of $\det( J_{g_{\bar\k}}(t))$ is not identically zero as long as  $\bar{\k}\neq 0$.
As the hypotheses of Theorem~\ref{thm:Kac-Rice} hold, 
the expected number of positive solutions to the system for parameters in $B$ is given by the Kac-Ric integral \eqref{eq:kac-rice-formula}
\begin{multline*} 
\hat{r}(B) = \int_0^{+\infty}\int_0^{+\infty}\int_{a_2}^{b_2}\int_{a_3}^{b_3}\int_{a_4}^{b_4}
\big| \det\big(J_{g_{\bar{\k}}}(t)\big) \big|\chi_{[a_1,b_1]}(g_{\bar\k,1}(t))\cdot\chi_{[a_5,b_5]}(g_{\bar\k,2}(t))\cdot \\ \ \left( \prod_{i=1}^5 \tfrac{1}{b_i-a_i}\right)d\k_2\, d\k_3\, d\k_4\, dt_2 \, dt_1.
\end{multline*}

	\begin{table}[t!]
			{\resizebox{\linewidth}{!}{
			\begin{tabular}{|c|c|c|c|c|c|c||c|c|c|c|c|c|}
				\hline
				\multirow{3}{*}{$N$} & \multicolumn{6}{c||}{(a) Uniform distribution} & \multicolumn{6}{c|}{(b)  Truncated normal distribution}  \\ \cline{2-13}
				& \multicolumn{3}{c|}{Simple Monte Carlo} & \multicolumn{3}{c||}{Antithetic Monte Carlo} & \multicolumn{3}{c|}{Simple Monte Carlo} & \multicolumn{3}{c|}{Antithetic Monte Carlo} \\ \cline{2-13}
				& $\widehat{I}_N$      & $\hat{e}_N$      &  Time     & $\widehat{I}_N$      &  $\hat{e}_N$     & Time  & $\widehat{I}_N$      & $\hat{e}_N$      &  Time   & $\widehat{I}_N$      &  $\hat{e}_N$     & Time   \\ \hline
				$10$ &  0.527     & 0.583      & 0.000023      & 0.199      & 0.172      & 0.000008      &  0.000     & 0.000      & 0.000016      & 0      & 0      & 0.000023      \\ \hline
				$10^2$ & 2.542      & 1.884      & 0.000021      & 1.100      & 0.415      & 0.000015    & 0.000      & 0.000      & 0.000060      & 0.000      & 0.000      & 0.000057        \\ \hline
				$10^3$ & 2.470      & 0.745      & 0.000165      & 0.942      & 0.147      & 0.000095    & 0.000      & 0.000      & 0.000543      & 0.000      & 0.000      & 0.000391        \\ \hline
				$10^4$ & 1.468      & 0.175      & 0.001570      & 1.662      & 0.333      & 0.000950    & 1.102      & 1.101      & 0.005587      & 0.002      & 0.002      & 0.004718        \\ \hline
				$10^5$ & 1.990      & 0.595      & 0.015469      & 1.392      & 0.055      & 0.008920    & 0.127      & 0.069      & 0.054673      & 2.027      & 0.817      & 0.038053        \\ \hline
				$10^6$ & 1.432      & 0.031      & 0.150625      & 1.449      & 0.034      & 0.090292    & 1.019      & 0.178      & 0.520236      & 1.021      & 0.171      & 0.372723       \\ \hline
				$10^7$ & \cellcolor[HTML]{FCFF2F}1.422      & \cellcolor[HTML]{FCFF2F}0.007      & \cellcolor[HTML]{FCFF2F}1.536203      & \cellcolor[HTML]{FCFF2F}1.419      & \cellcolor[HTML]{FCFF2F}0.009      & \cellcolor[HTML]{FCFF2F}0.939078   & 0.963      & 0.056      & 5.218110      & 0.965      & 0.057      & 3.681651         \\ \hline
				$10^8$ & 1.413      & 0.003      & 15.49415      & 1.415      & 0.003      & 10.04099    & 1.020      & 0.019      & 52.14723      & 1.034      & 0.019      & 37.95210        \\ \hline
				$10^9$ & 1.419      & 0.001      & 155.5443      & 1.418      & 0.001      & 92.33368    & \cellcolor[HTML]{FCFF2F}1.010      & \cellcolor[HTML]{FCFF2F}0.006      & \cellcolor[HTML]{FCFF2F}537.2604      & \cellcolor[HTML]{FCFF2F}0.989      & \cellcolor[HTML]{FCFF2F}0.006      & \cellcolor[HTML]{FCFF2F}371.1164        \\ \hline
		\end{tabular}}
		
		\medskip
			\caption{\small Output for network~\eqref{Equation:Joshi_2v5p_network} with computations done in \texttt{Julia}.  The time is reported in seconds and rounded. In (b) the minimum considered sample size is $10^5$ (see Subsection~\ref{sec:Numerical_integration}). The first cells satisfying the stop condition are highlighted. With two digits of significance, $\hat{r}$ is $1.4$ and $1.0$ in (a) and (b) respectively.
			}
			\label{Table:Efficiency_of_Antithetic}
			}
	\end{table}

As $0<t_1,t_2$  in $A$ and $t_1+t_2=\k_5$, the values of $t_1$ and $t_2$ as solutions to the system are bounded above by $b_5$.
In the computation of $\widehat{I}_N$ using Monte Carlo, we consider $\mu(t_1,t_2)$ to be the density of $U(0,b_5)\times U(0,b_5)$. 
Given sampled points $t_1^{(i)},t_2^{(i)},\k_2^{(i)},\k_3^{(i)},\k_4^{(i)}$  for $i=1,\dots,N$ and $N$ large, the Kac-Rice integral is approximated by the following sum:
\begin{equation}\label{eq:example_KR}
	\tfrac{b_5^2}{(b_1-a_1)(b_5-a_5)}\sum_{i=1}^N  \big| \det\big(J_{g_{\bar{\k}^{(i)}}}(t^{(i)})\big) \big| \chi_{[a_1,b_1]}\big(g_{\bar{\k}^{(i)},1}(t^{(i)})\big)\chi_{[a_5,b_5]}\big(g_{\bar{\k}^{(i)},2}(t^{(i)})\big).
\end{equation}
	
To illustrate this, consider the bounded box
\begin{equation}\label{eq:Joshi2v5p_box}
B= (0,100) \times  (0,2)\times  (0,200)\times  (0,100)\times  (0,2).
\end{equation}
Table~\ref{Table:Efficiency_of_Antithetic}(a) summarises the computed approximation of $\hat{r}(B)$   using \eqref{eq:example_KR} with simple and antithetic Monte Carlo, as $N$ is increased. 
This shows that the expected number of positive solutions to the system for parameters in $B$ is around $1.4$. 
Antithetic Monte Carlo is about $68$\% faster than simple Monte Carlo in this case with the same accuracy.

	\smallskip
	For comparison, we considered also truncated normal distributions on the parameters $\k_i\sim\bar{N}_{(a_i,b_i)}(\mu_i,\sigma_i)$ with  probability density function $\rho_i$, where the mean $\mu_i$ is the center of the interval $(a_i,b_i)$, and $\sigma_i=0.1$. Then the Kac-Rice integral  can be approximated with the following Monte Carlo sum
	\[b_5^2\sum_{i=1}^N  \big| \det\big(J_{g_{\bar{\k}^{(i)}}}(t^{(i)})\big) \big| \rho_1(g_{\bar\k^{(i)},1}(t^{(i)}))\rho_5(g_{\bar\k^{(i)},2}(t^{(i)})),\]
after sampling using $t_1\sim U(0,b_5),\;t_2\sim U(0,b_5)$ and $\k_i\sim \bar{N}_{(a_i,b_i)}(\mu_i,\sigma_i)$ for $i=2,3,4$. 
Results are shown in Table~\ref{Table:Efficiency_of_Antithetic}(b) for the box in \eqref{eq:Joshi2v5p_box}. Antithetic Monte Carlo is about $49$\% faster than simple Monte Carlo   in this case with the same accuracy.

\medskip
We have also used this example to compare the time it takes to compute $\widehat{I}_N$ and $\hat{e}_N$ using the algorithm \eqref{Equation:Standard_error_using_delta} on different platforms. In Table~\ref{Table:Comparison_table} we report the largest 
value of the type $N=10^d$ that can be computed under 200 seconds. Among the five considered platforms, \texttt{Numba} is the fastest.

\begin{table}[t!]
\renewcommand{\arraystretch}{1.2}
	{\small \begin{tabular}{c|c|c|c|}
		\cline{2-4}
		& \multirow{2}{*}{$N$} & \multicolumn{2}{c|}{Time} \\ \cline{3-4} &&
   Monte Carlo    & Antithetic Monte Carlo  \\ \hline
		\multicolumn{1}{|l|}{\texttt{Maple 2020}} & $10^5$ & 177.116 & 96.411 \\ \hline
		\multicolumn{1}{|l|}{\texttt{Python 3.7.4}} & $10^7$ & 89.74713 & 78.52874 \\ \hline
		\multicolumn{1}{|l|}{\texttt{C++11 Dev-Cpp 5.11}}   & $10^8$ & 142.395 & 94.6376   \\ \hline
		\multicolumn{1}{|l|}{\texttt{Numba 0.48.0}}   & $10^9$ & 93.41889 & 59.82621 \\ \hline
		\multicolumn{1}{|l|}{\texttt{Julia 1.4.2}}   & $10^9$ & 156.2837 & 88.81576 \\ \hline
		\multicolumn{1}{|l|}{\begin{minipage}[h]{0.28\textwidth}\vspace{0.1cm}   \texttt{Julia 1.4.2} parallelized \\ with 2 workers  \vspace{0.1cm} \end{minipage}} &  $10^9$ &  76.16551 & 47.35521 \\ \hline
		\multicolumn{1}{|l|}{\begin{minipage}[h]{0.28\textwidth}\vspace{0.1cm}   \texttt{Julia 1.4.2} parallelized \\ with 32 workers  \vspace{0.1cm} \end{minipage} } &  $10^{10}$ & 116.7884 & 34.55502 \\
\hline
	\end{tabular}
	}

	\medskip
	\caption{\small Computation time  (in seconds) for $\widehat{I}_N$ in \eqref{eq:example_KR} and $\hat{e}_N$ using \eqref{Equation:Standard_error_using_delta}. Computation is performed using both simple and antithetic Monte Carlo sampling and different platforms. 
	We report the computation time  of the largest sample size $N=10^d$ taking less than 200 seconds. The antithetic Monte Carlo is faster on all  platforms,  and \texttt{Numba} is the  fastest of the considered platforms. The parallel implementation in \texttt{Julia} using the package \texttt{Distributed} with $d$ workers increases the speed of computations by a factor of $\max(d,\widetilde{d}/2)$ where $\widetilde{d}$} is the number of cpus of the computer.
	\label{Table:Comparison_table}
\end{table}

\medskip
\paragraph{\bf Finding a point in multistationary region.}
This network admits between one and three positive steady states. 
We use the Kac-Rice integral and Monte Carlo integration (in \texttt{Julia}) to find a parameter point where the network has three positive steady states. To this end, we consider the following starting box:
\[ B= (0,100)\times (0,2) \times (0, 200) \times (0,100) \times (0,2).\]
The algorithm outlined for \textbf{Problem II} concludes with the box
\[ C=(50,75)\times (1,1.5)\times (150,200)\times (0,50)\times (1,2),\]
after $8$ steps in less than $7.5$ seconds. 
Therefore, for almost all parameters in $C$, the network has three positive steady states, and hence multistationarity.

 \subsection{Finding $5$ solutions}\label{Example:Degree5}
 We now analyse an example where the maximal number of solutions is five and $A$ is unbounded. 
 Consider the following parametrized univariate polynomial of degree five in the variable $t$ and parameters $\k_1,\k_2$:
	\begin{equation}\label{Equation:f_Degree5}
	\begin{split}
	f_{\k}(t) &= t^5-(\k_1+\tfrac{9}{2})t^4+(\tfrac{9}{2}\k_1+\tfrac{21}{4})t^3+(-\tfrac{23}{4}\k_1+\tfrac{3}{8})t^2\\ & \hspace{5cm}+(\tfrac{15}{8}\k_1-\tfrac{23}{8})t+(\tfrac{1}{100}\k_2-\tfrac{1}{16}).
	\end{split}
	\end{equation}
	
\begin{figure}[t!]
		\centering
		\hspace{-1cm}
		\begin{minipage}[h]{0.29\textwidth}
		\centering
			\includegraphics[width=\linewidth]{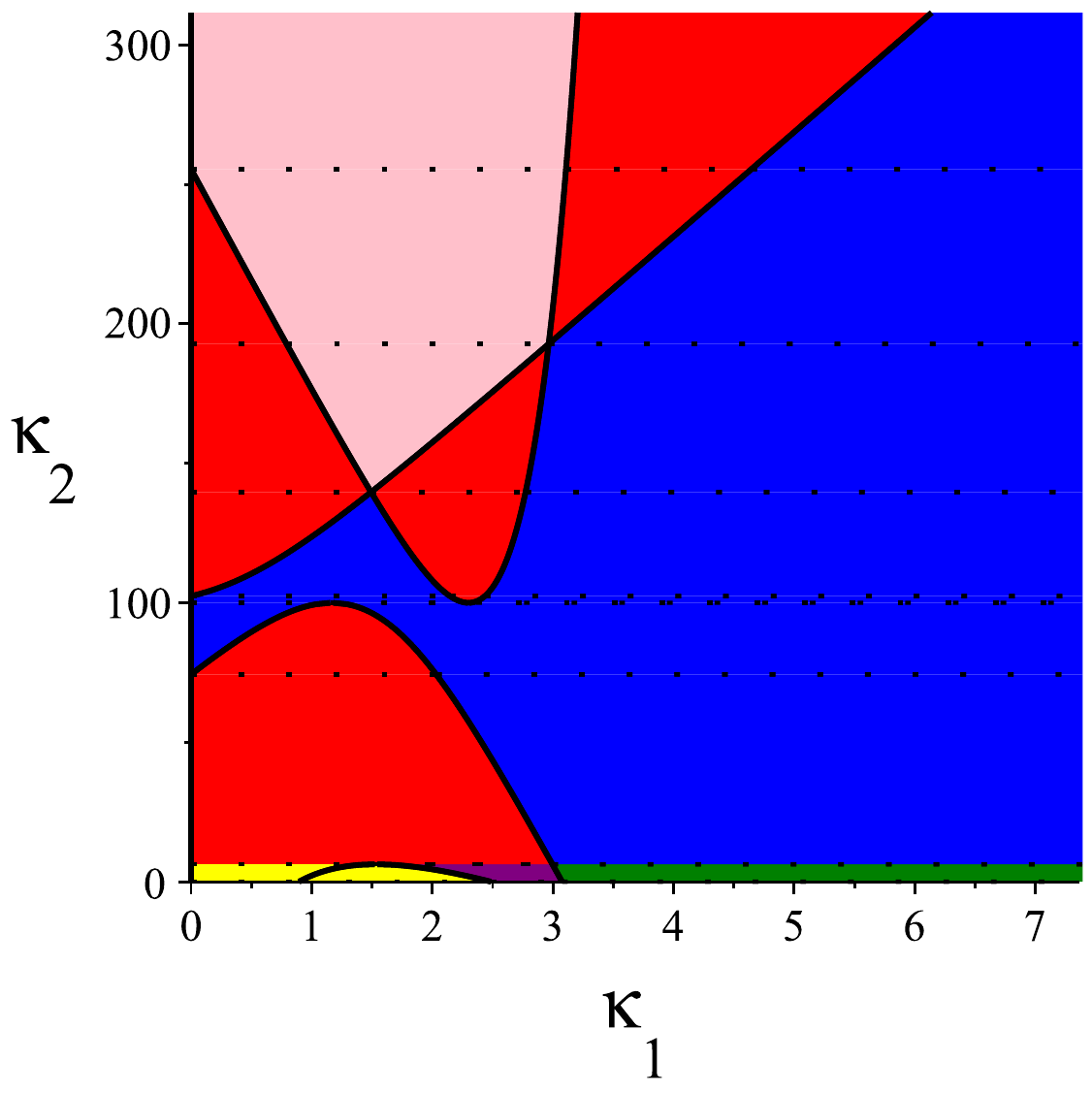}
			
			(a)
		\end{minipage}
		\begin{minipage}[h]{0.29\textwidth}
		
				\centering
			\includegraphics[height=\linewidth]{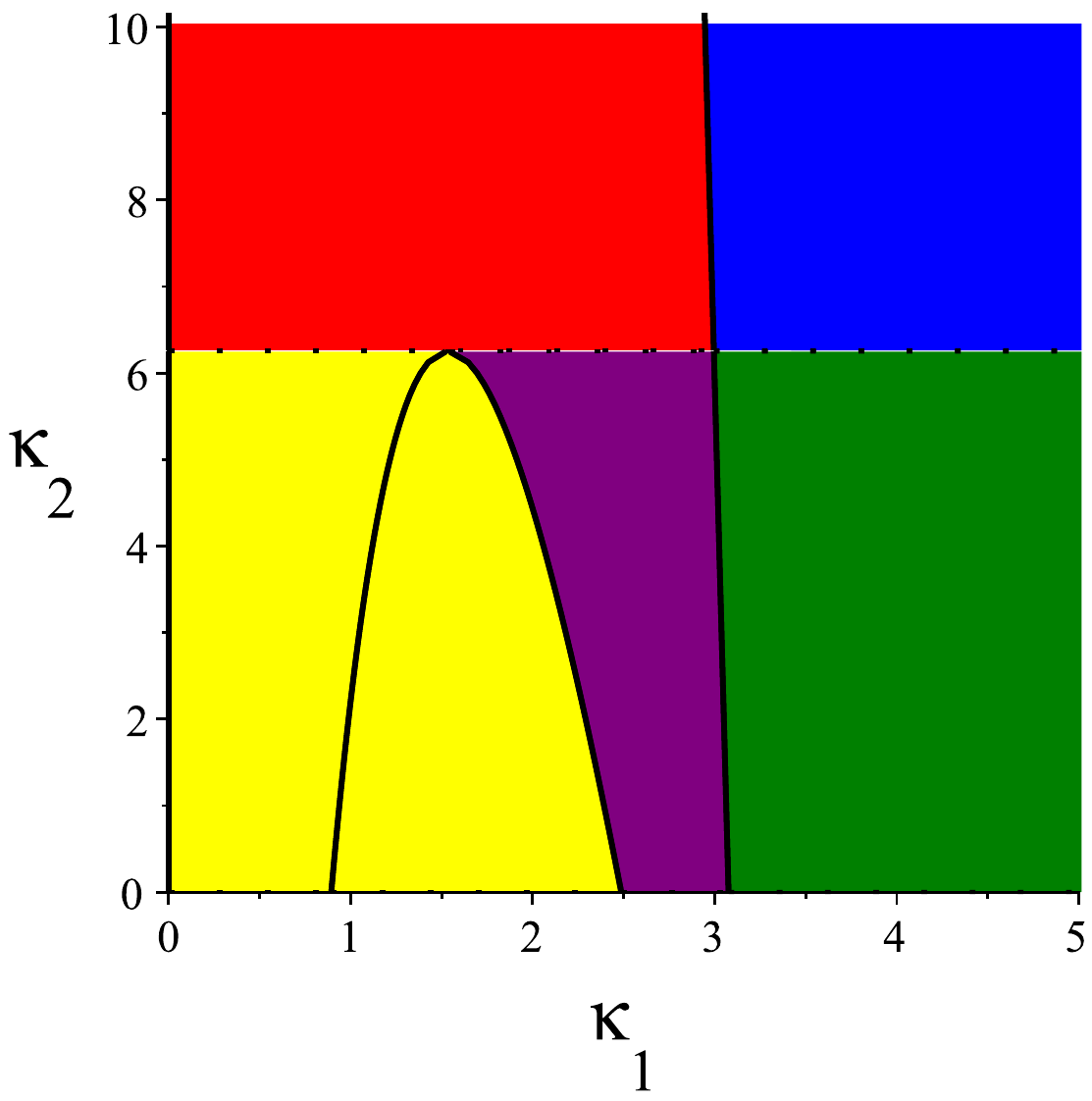}
			
			(b)
		\end{minipage}
		\begin{minipage}[h]{0.29\textwidth}
		
				\centering
			\includegraphics[height=\linewidth]{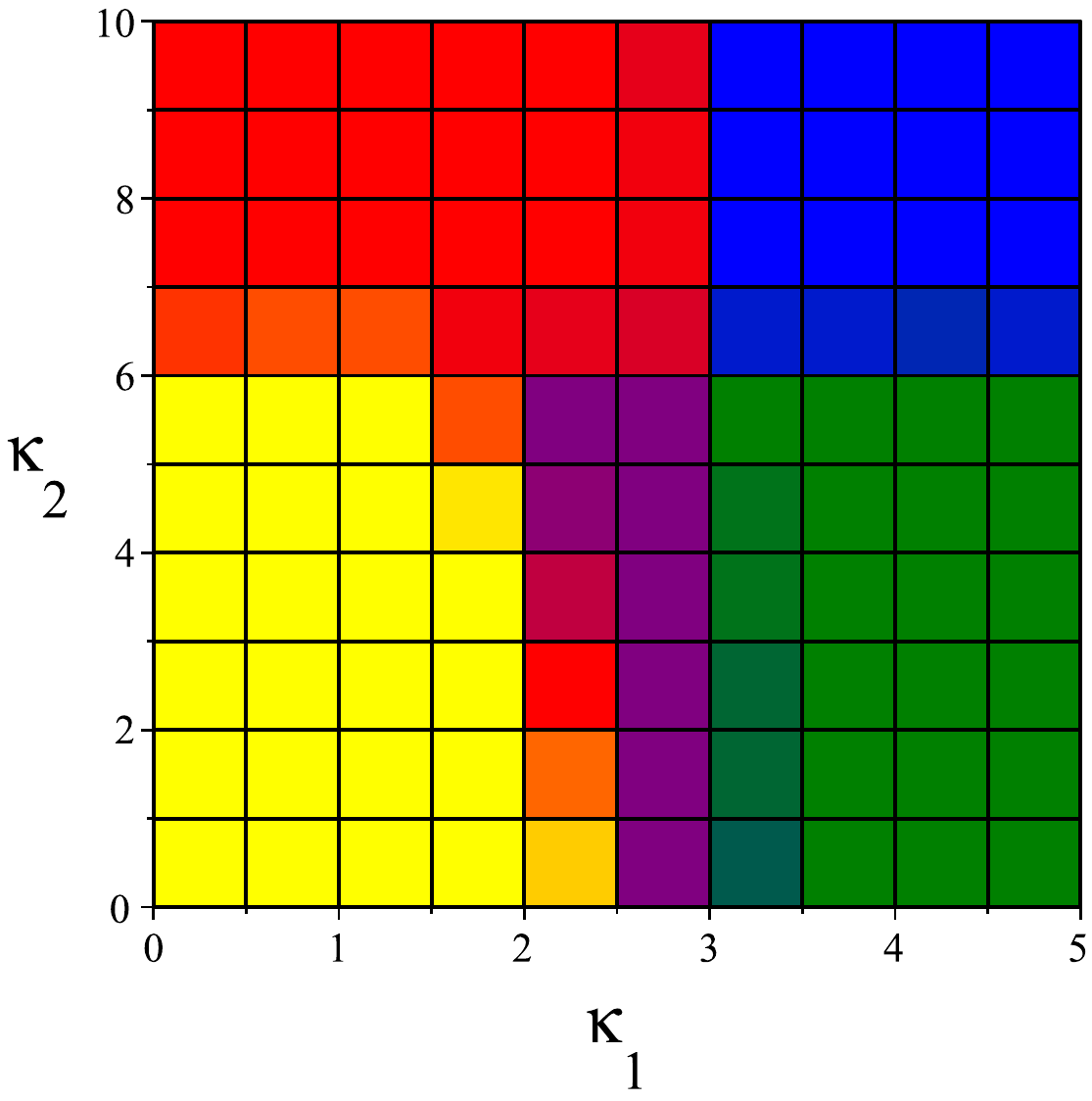}
			
			(c)
		\end{minipage}
		\begin{minipage}[h]{0.06\textwidth}
			
			\centering
			\vspace{-1cm}\includegraphics[height=3.5cm]{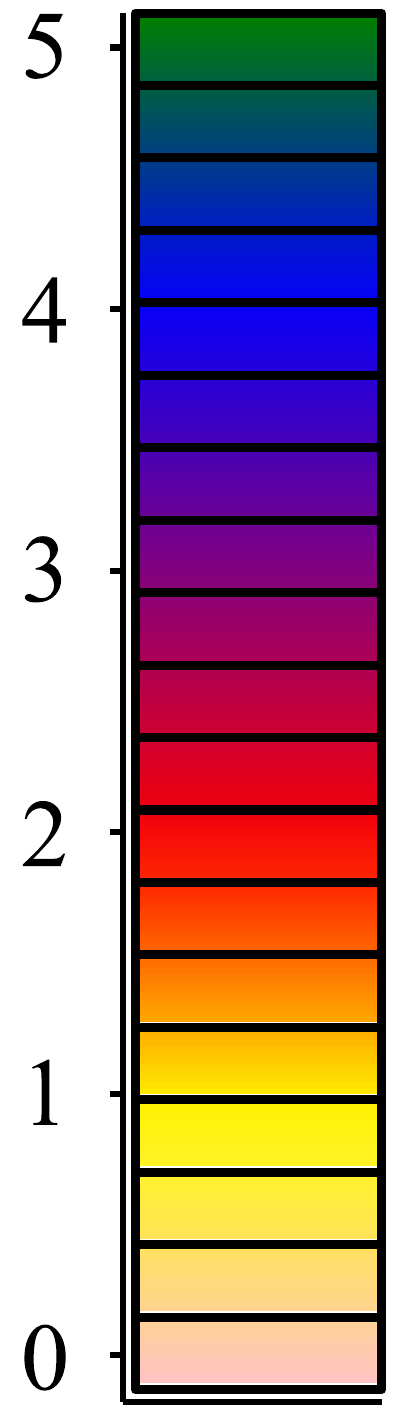}
			
		\end{minipage}
		\caption{{\small 
		Partition of the parameter region of \eqref{Equation:f_Degree5} according to the number of positive roots. The color indicates the number of positive roots as given by the bar code on the right. 
		(a-b) Partition obtained using CAD for two different boxes.  (c) Grid partition and computation of $\hat{r}$ with the Kac-Rice formula and Monte Carlo integration.}} 
		\label{fig:degree5_CAD_all}\label{fig:degree5_CAD}\label{fig:Example4}\label{Figure:Degree5}
	\end{figure}

	Using  CAD, we know that $f_{\k}(t)$ generically  has 0, 1, 2, 3, 4 or 5 positive roots for suitable choices of the parameter vector $(\k_1,\k_2)\in\R_{>0}^2$ (see Figure~\ref{fig:degree5_CAD_all}(a-b)).  The polynomial $f_{\k}(t)$ is linear in $\k_2$ with coefficient $h(\k_1,t)=100$, and the hypotheses of Theorem~\ref{thm:Kac-Rice} hold for $A=\R_{>0}$. Hence, for a box $B=[a_1,b_1]\times [a_2,b_2]$ with $0\leq a_i$, $\hat{r}(B)$ is given by:
	\begin{equation}\label{Equation:Degree5_1}
	\hat{r}(B)=\tfrac{1}{b_2-a_2} \int_0^{+\infty}\int_{a_1}^{b_1}\big|g'_{\k_1}(t)\big|\chi_{[a_2,b_2]}(g_{\k_1}(t))\,d\k_1\, dt,
	\end{equation}
	where
	\begin{equation*}
	g_{\k_1}(t)=-100\big(t^5+(-\k_1-\tfrac{9}{2})t^4+(\tfrac{9}{2}\k_1+\tfrac{21}{4})t^3+(-\tfrac{23}{4}\k_1+\tfrac{3}{8})t^2+(\tfrac{15}{8}\k_1-\tfrac{23}{8})t-\tfrac{1}{16}\big).
	\end{equation*}

In order to approximately partition the box $B=[0,5]\times [0,10]$ according to the number of positive roots of $f_\k(t)$, we subdivide it into 100 sub-boxes and compute $\hat{r}$ using Monte Carlo integration.
 By using the transformation
\[
	\int_0^{+\infty} h(t)dt = \int_0^1 h(t)dt+\int_0^1\tfrac{1}{t^2}h(\tfrac{1}{t})dt
\]
where $h(t)$ is the integrand of  \eqref{Equation:Degree5_1} in $t$, the integral  \eqref{Equation:Degree5_1} 
can be written as the following  integral:
	\begin{equation*}
	\hat{r}(B)=\tfrac{1}{b_2-a_2} \int_0^1\int_{a_1}^{b_1}\Big(\big|g'_{\k_1}(t)\big|\chi_{[a_2,b_2]}(g_{\k_1}(t))+\tfrac{1}{t^2}\big|g'_{\k_1}(\tfrac{1}{t})\big|\chi_{[a_2,b_2]}(g_{\k_1}(\tfrac{1}{t}))\Big)\,d\k_1\, dt.
	\end{equation*}
	By choosing $t$ and $\k_1$ to follow uniform distributions on $[0,1]$ and $[a_1,b_1]$ respectively, $\hat{r}(B)$ is approximated by 		\begin{equation*}
	\tfrac{1}{b_2-a_2}\sum_{i=1}^N \Big(\big|g'_{\k_1^{(i)}}(t^{(i)})\big|\chi_{[a_2,b_2]}(g_{\k_1^{(i)}}(t^{(i)}))
+\tfrac{1}{(t^{(i)})^2}\big|g'_{\k_1^{(i)}}(\tfrac{1}{t^{(i)}})\big|\chi_{[a_2,b_2]}(g_{\k_1^{(i)}}(\tfrac{1}{t^{(i)}}))\Big)
	\end{equation*}
for  sampled points $t^{(i)}$, $\k_1^{(i)}$  for $i=1,\dots,N$ and $N$ large.

	Figure~\ref{Figure:Degree5}(c) depicts each of these sub-boxes, colored corresponding to the  approximated value of $\hat{r}$ for $N=10^{9}$. Note that $\hat{r}(C)\simeq 2$ for the box $C=[2,2.5]\times [2,2.5]$. However this sub-box is not located inside or even have intersection with the (open) parameter region where $f_{\k}(t)$ has two positive roots. It intersects only regions with one and three solutions, but the areas of the two intersections are almost equal. 
Only when $\hat{r}(C)$ is zero or five ($M_{\rm min}$ and $M_{\rm max}$ here), we can conclude that almost all parameter points in the box yield to zero or five positive roots. For example, the box	$[3.5,5]\times [0,6]$ is entirely inside the parameter region with five positive roots. 

For this computation,  the standard error $\hat{e}_N$ increased
with $a_1$, going from $0.002$ to a maximal value of $0.02$ for the box $[4.5,5] \times [9,10]$ independently of $a_2$.

\subsection{Dual phosphorylation }
\label{sec:twosite} 
We consider   the following reaction network:
\begin{align*}
	X_1+X_4\ce{<=>[k_1][k_2]}X_6\ce{->[k_3]}X_2 &+ X_4\ce{<=>[k_4][k_5]}X_7\ce{->[k_6]}X_3+X_4\\
	X_3+X_5\ce{<=>[k_7][k_8]}X_8\ce{->[k_9]}X_2 &+ X_5\ce{<=>[k_{10}][k_{11}]}X_9\ce{->[k_{12}]}X_1+X_2.
\end{align*}
This network models the distributive and sequential phosphorylation and dephosphorylation of a substrate that is phosphorylated at none, one or two sites  ($X_1$, $X_2$ and $X_3$), catalysed by enzymes $X_4,X_5$.  Generically, this network has either one or three positive steady states \cite{Wang:2008dc}. 
Even though substantial work has been done to understand the parameter region of multistationarity \cite{maya-bistab,dickenstein:regions,feliu-twosite,bates-gunawardena,Conradi_Iosif_Kahle}, an explicit description is still unknown.

This network has three conservation laws
\begin{align*}
	x_1+x_2+x_3+x_6+x_7+x_8+x_9 &= T_1, &
	x_4+x_6+x_7 &= T_2, &
	x_5+x_8+x_9 &= T_3.
\end{align*}
System \eqref{eq:steadystates} can be simplified to a system of three polynomials in the three variables $x_1,x_4,x_5$ and $15$ parameters, see \cite{dickenstein:regions}. 
Essentially, the equations $\widetilde{F}_{k}(x)=0$ in \eqref{eq:steadystates} can be solved for 
$x_1,x_4,x_5$, and the output inserted into the three conservation laws. Hence the three polynomials 
are  linear in $T_1,T_2,T_3$, respectively. With $\bar{\k}=(k_1,\dots,k_{12})$,  the hypotheses of Theorem~\ref{thm:Kac-Rice} hold. 
The numerator of the rational function $\det(J_{g_{\bar{\k}}}(t))$ has total degree $18$ (in the variables and parameters) and   $165$ terms. The denominator has total degree $10$ and $9$ terms. 

By \cite[Corollary~4.13]{Conradi_Iosif_Kahle}, parameter choices where $T_1<T_2$ and $T_1<T_3$ do not yield  multistationarity. Therefore for the following box  $\hat{r}(B)$  is one:
\begin{multline*}
B = (500,1000) \times (25,50) \times (25,50) \times (5,10) \times (5,10) \times (5,10) \times (5,10) \times (1,2)\\
\times (1,2) \times (5,10) \times (50,100) \times (50,100) \times (1,2) \times (2,4) \times (2,4).
\end{multline*}
With Monte Carlo integration, it takes 3.2 seconds to approximate $\hat{r}(B)$ to  one with two digits of significance  using $N=10^7$ (the minimum sample size is set to $10^3$). 

We consider  this other box:
\begin{multline*}
B = (0.5, 1.5) \times (509.5, 510.5) \times (1.5, 2.5) \times (1.5, 2.5) \times (0.5, 1.5) \times (0.5, 1.5)\\ \times 
(1.5, 2.5) \times (0.5, 1.5) \times (0.5, 1.5) \times (1.5, 2.5) \times (0.5, 1.5) \times (0.5, 1.5) \times (110, 150)\\
\times (20, 30) \times (15, 25).
\end{multline*}
In this case, the minimum sample size for the computation of $\hat{r}(B)$ is $N=10^7$.
With antithetic Monte Carlo we obtain that $\hat{r}(B)=1.45$ in $27564$ seconds   with standard error $\hat{e}=0.015$ (1011.5 seconds using $32$ workers). Therefore this box intersects the  region of multistationarity.

In order to find a box inside the region of multistationarity, we have considered the bisect approach with the box $B$ and parallelized with $32$ workers. After $25$ bisect steps, the computation of $51$ integrals, and approximately $16$ hours (precisely $58{,}156.124$ seconds), we obtain the box
\begin{multline*}
(1.125, 1.25) \times(510.0, 510.25) \times(1.5, 1.75) \times(2.25, 2.5)   \times(0.5, 0.75) \times(0.5, 0.75) \times(2.25, 2.5)  \\ \times(0.5, 0.75) \times(0.5, 0.75) \times(2.25, 2.5) \times(0.5, 1.0) \times(1.0, 1.5) \times(115.0, 120.0) \\  \times(25.0, 30.0) \times(20.0, 25.0).
\end{multline*}
At the end of the computation, this box gives $\hat{r}$ equal to $3.01$ and $\hat{e}=0.05$, with $N=10^{11}$.
An additional computation with a higher sample size $N=10^{12}$, gives $\hat{r}= 2.94$ with $\hat{e}= 0.017$ (in $11393$ seconds).
This implies that a big portion of this box is inside the multistationarity region. 
 Note that numerical error implies that different runs of the bisect strategy might yield to different boxes of multistationarity.

\subsection{Extended hybrid histidine-kinase network }\label{sec:extendedHK}
Finally, as a last example, we consider an extension of the hybrid histidine-kinase network studied in Subsection~\ref{Example:HK_network} as given in \cite{FeliuHK}:
\begin{equation}\label{eq:extendedHK}
\begin{aligned}
X_1\ce{->[k_1]}X_2 & \ce{->[k_2]}X_3\ce{->[k_3]}X_4 & \qquad X_7\ce{->[k_7]}X_8 & \ce{->[k_8]}X_9\ce{->[k_9]}X_{10}\\
X_3+X_5 & \ce{->[k_4]}X_1+X_6 & X_9+X_5 & \ce{->[k_{10}]}X_7+X_6\\
X_4+X_5 & \ce{->[k_5]}X_2+X_6 & X_{10}+X_5& \ce{->[k_{11}]}X_8+X_6 \\
	X_6& \ce{->[k_6]}X_5.
\end{aligned}
\end{equation}
This network has three conservation laws
\begin{align*}
	x_1+x_2+x_3+x_4 &= T_1, &
	x_5+x_6 &= T_2, & 
	x_7+x_8+x_9+x_{10} &= T_3,
\end{align*}
and hence  involves a $14$-dimensional parameter vector: $\k=(k_1,\dots,k_{11},T_1,T_2,T_3)$. 
By \cite{FeliuHK},  network~\eqref{eq:extendedHK} admits between one and five positive steady states. The corresponding system \eqref{eq:steadystates} can be simplified to a univariate polynomial $f_\k(t)$  of degree five in $t=x_5$, which further is linear in each of $\k_{12}$, $\k_{13}$ and $\k_{14}$. 
We choose to isolate $\k_{12}=T_1$, such that $\bar{\k}=(k_1,\dots,k_{11},T_2,T_3)$. The hypotheses of Theorem~\ref{thm:Kac-Rice} hold.

We use Monte Carlo integration and the Kac Rice formula to address \textbf{Problem~II} and find a multistationary point in the box
\begin{multline}\label{eq:box_doubleHK}
B=(0,2)\times (0,100)\times (0,50)\times (0,10)\times (0,100)\times (0,50)\times (0,5)\\ \times (50,100)\times (0,100)
\times (0,100)\times (0,100)\times (0,100)\times (0,50)\times (0,100).
\end{multline}
We have $\hat{r}(B)=1.21$ with $\hat{e}=0.010$.
We set the termination condition of the search algorithm to be $\hat{r}(C)\geq 2.95$. After $32$ bisections, computing $65$ integrals in $3083$ seconds, the algorithm terminates and returns the following sub-box:
\begin{multline*}
C=(0.25,0.5)\times (75,87.5)\times (43.75,50)\times (3.75,5)\times (87.5,100)\times (18.75,25)
\times (1.25,2.5) \\\times (62.5,75.0)\times (75,100)\times (25,50)\times (75,100)\times (75,100)
\times (25,50)\times (50,100).
\end{multline*}
As  $\hat{r}(C)=2.95$ with $\hat{e}_N=0.008$, 
there must be parameter points yielding to three or five positive steady states.
For $975$ out of $1000$ random parameter points in $C$, the polynomial has three positive roots (found numerically).

\begin{figure}[t]
	\centering
	\includegraphics[width=4cm]{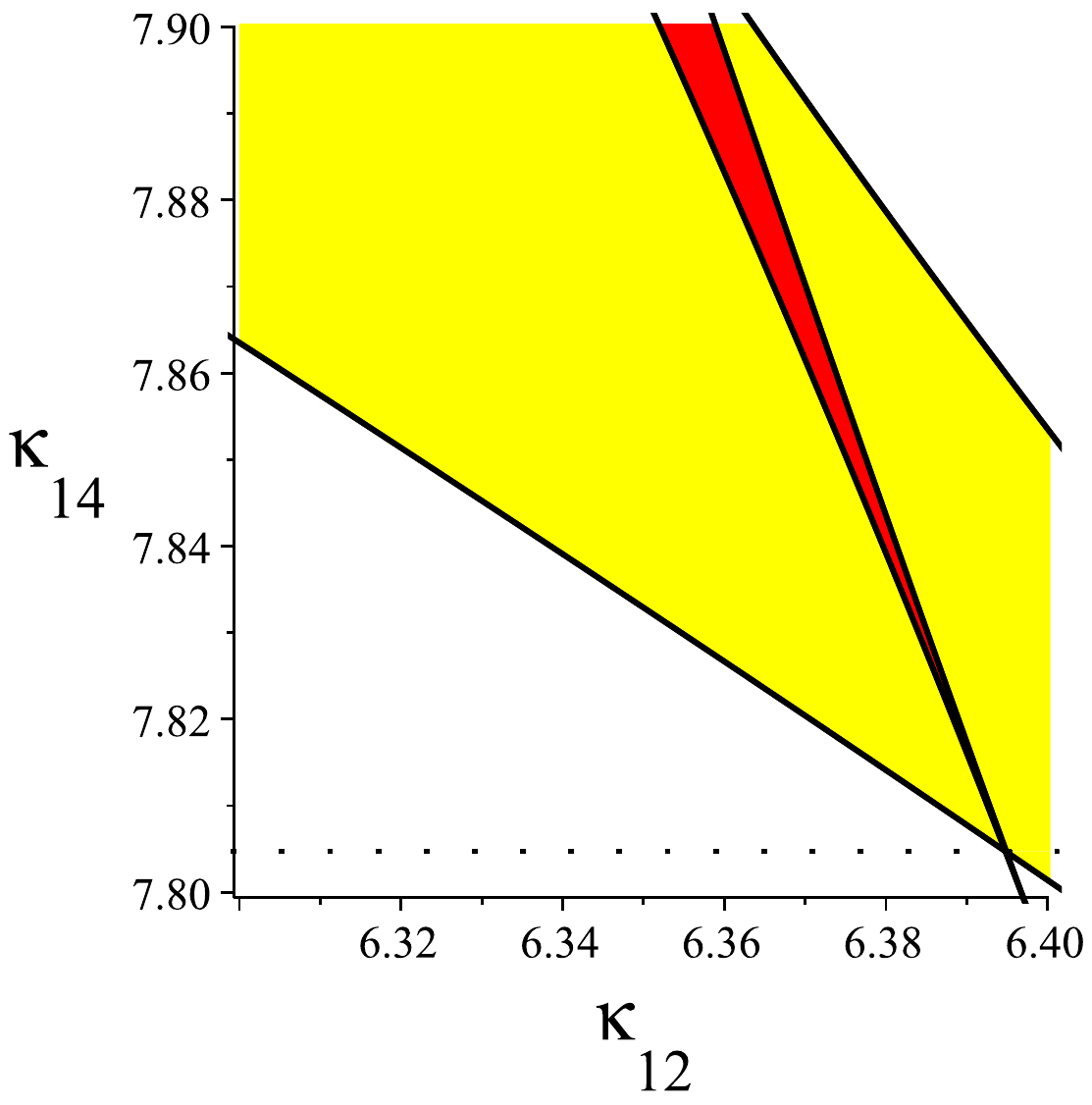}
	\caption{\small In the red, yellow and white regions, network~\eqref{eq:extendedHK} with the choices in  \eqref{Equation:DoubleHK_parameters} has five, three and one positive steady states respectively, as obtained using CAD.} 
	\label{Figure:DoubleHK_CAD}
\end{figure}

For the box $B$,  the minimum sample size for Monte Carlo integration is $10$. However, it is not always the case. For example,  fix all parameters other than $\k_{12}$ and $\k_{14}$ as follows 
\begin{equation}
\begin{aligned}
	\k_1 &= 0.1, & \k_2 &= 120, & \k_3 &= 17.95, & \k_4 &= 0.1795, & \k_5 &= 0.713, & \k_6 &= 1,\label{Equation:DoubleHK_parameters}\\
	\k_7 &= 0.002, & \k_8 &= 500, & \k_9 &= 160, & \k_{10} &= 0.147, & \k_{11} &= 4.15, & \k_{13} &=16.27.
	\end{aligned}
\end{equation}
These values are taken from \cite[Fig.~3B for $n=2$]{FeliuHK}, and, for the right choice of $\k_{12}, \k_{14}$, yield to five positive steady states. 
We have
	{\small \begin{align*}
g_{\k_{14}}(t)& = 93770052422884376700t^5   + (18753935468835 \k_{14}     - 146395097463035713.8909 )10^4 t^4 \\ &   + (204244474710835\k_{14}   +360657036215560.9291) 10^6t^3\\& 
	 + ( 1332479842\k_{14} - 221423667689.22618)10^{11}t^2 \\ & +  (2860512\k_{14} -206772.28792)10^{15}t -112145856\cdot 10^{14}\; /\\
	&\qquad \big(468.4598789t^4 +46854.40101t^3 +1087.040656t^2 + 24572.832 t\big).
\end{align*}}%
The coefficients of $g_{\k_{14}}(t)$ are of different scales ranging from $10^{23}$ to $10^3$. We use Monte Carlo integration to approximate the average number of positive steady states when $(\k_{12},\k_{14})$ belong to the following box:
\[B=[6.3,6.4]\times [7.8,7.9].\]
For any sample size from 10 to $10^{11}$, we obtain $\widehat{I}_N=\hat{e}_N=0$. It is clear from Figure~\ref{Figure:DoubleHK_CAD} that   $1<\hat{r}(B)$, and hence different from zero. 
At this point it is unclear to us whether the problem arises because the minimum sample size is larger than $10^{11}$, or due to numerical errors arising from the different orders of the coefficients of $g_{\k_{14}}$ and limited machine number sizes.

For comparison,  instead of the values in \eqref{Equation:DoubleHK_parameters}, we let all the parameters other than $\k_{12}$ and $\k_{14}$ be equal to 1. Then for the box $B=[0,1]\times [0,1]$, we obtain $\widehat{I}_N=1$ with $N=10^5$ in $0.05$ seconds and two digits of precision. Here $N=10$ is already acceptable as  sample size, and the function  $g_{\k_{14}}(t)$ does not have coefficients of different order:
\begin{align*}
g_{\k_{14}}(t) &= \tfrac{-4t^5 + (2\k_{14} + 4)t^4 - 4\k_{14}t^3 + (3\k_{14} - 4)t^2 + (\k_{14} - 3)t + 1}{2t^4 + 4t^3 + 3t^2 + t}. 
\end{align*}

\medskip
For this example, the parameter region with five steady states  is likely  too small to be detectable using our approach. For instance, for the boxes
$B=[0,100]^{14}$ or $[0.1,1]^{14}$, we obtain that $\hat{r}(B)=1$ with 3 digits of significance, meaning the regions with three or five steady states are small. 
However, we illustrated with the box in \eqref{eq:box_doubleHK} that parameter points yielding multistationarity can be easily found in a box of interest, even with $14$ free parameters.

\section{Proof of Theorem~\ref{thm:Kac-Rice}}\label{sec:proofKacRice}

In this section we prove  Theorem~\ref{thm:Kac-Rice}. The argument follows 
standard approaches to establish Kac-Rice formulas, see for example \cite{LevelSetsAndExtremaOfRandomProcessesAndFields}, Chapter 3 for the one-dimensional case, and Chapter 6 for the multivariate case.

Given a set $A\subseteq \R^n$, we let $A^\circ$, $\bar{A}$ and $\partial A$ denote respectively the interior, closure and boundary of $A$ (with respect to the Euclidean topology). 
Given a sequence of sets $\lbrace S_r\rbrace_{r\in\N}$ such that $S_1\subseteq S_2\subseteq S_3\subseteq \cdots$ and   $\cup_{r\in\N}S_r=S$, then   we use the notation $\lim_{r\rightarrow +\infty} S_r=\cup_{r\in\N}S_r=S$. 

Let $\pi_1$ and $\pi_2$ be respectively the projections of $\R^{m}$ onto the first $n$ components and the last $m-n$ components (so that $\bar{\k}=\pi_2(\k)$). 
For $t\in A$, it holds
\begin{equation}\label{eq:f_to_g}
f_\k(t)=0 \quad\textrm{ if and only if }\quad \pi_1(\k)= g_{\pi_2(\k)}(t).
\end{equation}

Before proving Theorem~\ref{thm:Kac-Rice}, we 
establish a series of lemmas. 

\begin{lem}\label{lem:refined}
With the notation and assumptions of Theorem~\ref{thm:Kac-Rice} 
the following holds:
\begin{enumerate}[(i)]
\item For every $\bar{\k}\in \widetilde{B}$ outside a Zariski closed set 
of measure zero $\widetilde{P}' \subseteq \widetilde{B}$, 
 there exists a measure zero set $A(\bar{\k})\subseteq A$ such that $\rho_i$ is continuous at $g_{\bar{\k},i}(t)$ for all $t\notin A(\bar{\k})$ and for all $i=1,\dots,n$.
\item With $\widetilde{P}\subseteq \widetilde{B}$ as in Theorem~\ref{thm:Kac-Rice}(ii), for all $\bar{\k}\notin \widetilde{P}$ there exists a Zariski closed set 
of measure zero $Q_{\bar{\k}}\subseteq \pi_1(B)$, such that 
if $u \in \pi_1(B) \setminus Q_{\bar{\k}}$, then 
the solution set to $g_{\bar{\k}}(t)=u$ in $A$ 
consists of a finite number of simple points in the interior of $A$. 
In particular, $\det(J_{g_{\bar{\k}}}(t))\neq 0$ at all solution points.
\end{enumerate}
\end{lem}
\begin{proof}
To show (i), we use that 
$\rho_i$ is a continuous function except maybe in a finite number of points $\xi_1,\dots,\xi_\ell$. 
For a fixed $j\in \{1,\dots,\ell\}$ and $\bar{\k}\in \widetilde{B}$, as $g_{\bar{\k}}$ is rational,  the solutions to $g_{\bar{\k}}(t)=\xi_j$ in $t$ form an algebraic variety given by $n$ equations and $n$ variables. As for almost all $\bar{\k}$,   $g_{\bar{\k}}$ is not constant (by assumption Theorem~\ref{thm:Kac-Rice}(ii)), for almost all $\bar{\k}$ this algebraic variety has codimension at least $1$ in $\mathbb{R}^n$. Hence (i) holds.

We turn now to (ii).  
By assumption Theorem~\ref{thm:Kac-Rice}(ii), for all $\bar{\k}\notin \widetilde{P}$, the polynomial $p_{\bar{\k}}(t)$ given by the numerator of $\det (J_{g_{\bar{\k}}}(t))$ is not identically zero. 
For a fixed  $\bar{\k} \notin \widetilde{P}$,  $p_{\bar{\k}}(t)=0$ defines a real algebraic variety $Y$ of codimension at least $1$ in $\mathbb{R}^n$, and hence is not Zariski dense in $A$. 
Define $Q_{\bar{\k}}$ as the image of $g_{\bar{\k}}$ restricted to $Y$, which is not Zariski dense in $\pi_1(B)$. 
For any $u$ outside  $Q_{\bar{\k}}$, any solution $t^*$ to $g_{\bar{\k}}(t)=u$  satisfies $p_{\bar{\k}}(t^*)\neq 0$, and hence is simple and isolated. In this case there is a finite number of solutions, as $g_{\bar\k}$ is rational.

As $A$ is a box, the boundary of $A$ can be decomposed into the union of subsets of coordinate hyperplanes. The restriction of $g_{\bar{\k}}$ to a coordinate hyperplane yields a rational function in  (at most) $n-1$ variables and $n$ entries, and the image is not Zariski dense in $\pi_1(B)$. 
Now augment $Q_{\bar{\k}}$ to include the image of $g_{\bar{\k}}$ restricted to every coordinate hyperplane describing the boundary of $A$. Then the equation $g_{\bar{\k}}(t)=u$ has no boundary solutions if $u\notin Q_{\bar{\k}}$. Finally, redefine $Q_{\bar{\k}}$ to be its Zariski closure, which by construction is a real algebraic variety different from $\R^n$ and hence has measure zero.
This concludes the proof of (ii).
\end{proof}

 \begin{lem}\label{Lemma_Kac_Rice_Unbounded_A}
 Let $f_\k(t)$, $B$,  and $A\subseteq \R^n$ as in Theorem~\ref{thm:Kac-Rice}. 
Assume that for any compact box $S\subseteq A$, the Kac-Rice formula \eqref{eq:kac-rice-formula} holds for the domain $S$.
Then the Kac-Rice formula \eqref{eq:kac-rice-formula} holds for $A$. 
\end{lem}
\begin{proof}
As $A=I_1\times \dots \times I_n$ with $I_j$ intervals, one can easily construct compact sets $S_{i,j}\subseteq I_j$ for $i\in \N$ such that $\lim_{i\rightarrow +\infty}S_{i,j} = I_j$ for all $j=1,\dots,n$. 
By definition the compact set $S_i:= S_{i,1}\times \dots \times S_{i,n}\subseteq A$,  and 
 $ \lim_{i\rightarrow +\infty}S_i=A$.

By hypothesis, formula \eqref{eq:kac-rice-formula} holds for the compact set $S_i$. 
Let $\psi(t)$ be the function integrated on the right-hand side of \eqref{eq:kac-rice-formula}  such that 
\begin{equation*}
\E\big(\#\big(f_\k^{-1}(0)\cap S_i\big)\big)=
\int_{S_i} \psi(t) dt \qquad \textrm{for all } i\in\N.
\end{equation*}
By the Lebesgue's monotone convergence theorem \cite[Theorem 1.26]{Rudin_Real_and_Complex}, we have 
\begin{equation*}
\int_{A} \psi(t) dt  = \lim_{i\rightarrow +\infty}  \int_{S_i} \psi(t) dt \qquad \textrm{and}\qquad 
\E\big(\#\big(f_\k^{-1}(u)\cap A\big)\big)=\lim_{i\rightarrow +\infty}\E\big(\#\big(f_\k^{-1}(u)\cap S_i\big)\big),
\end{equation*}
and hence
\[ \E\big(\#\big(f_\k^{-1}(u)\cap A\big)\big) = \int_{A} \psi(t) dt\]
as desired.
\end{proof}

 \medskip
Given $x\in \R^n$ and $\delta>0$, let $B(x,\delta)$ be the open ball centered at $x$ of radius $\delta$ and $V_\delta$  the volume of any such  ball of radius $\delta$.

\begin{lem}\label{Lemma:Dirac_to_Lebesgue3}
	Let $T\subseteq\R^n$ be an open set, $f\colon T\rightarrow\R$ a continuous function, and consider an increasing sequence of open sets $\{C_r\}_{r\in \N}$ such that $\bigcup_{r\in \N} C_r=T$.

Then, for $y\notin \partial T$, it holds
	\begin{equation*}
\lim_{r\rightarrow +\infty}\int_{C_r} \frac{ \chi_{B(x,1/r)}(y)}{V_{1/r}}f(x)\, dx = \begin{cases} f(y) &\textrm{if } y\in T, \\
0 &\textrm{ otherwise}.\end{cases}
	\end{equation*}
\end{lem}
\begin{proof}
If $y\notin T$, then as $y\notin \partial T$,  we have $y\in (\R^n\setminus T)^\circ$ and $y\notin B(x,\tfrac{1}{r})$ for all   $x\in T$ and $r$ large enough. 
Hence $\lim_{r\rightarrow +\infty} \int_{C_r}\frac{ \chi_{B(x,1/r)}(y)}{V_{{1/r}}}f(x)dx =0$.

Assume now $y\in T$. 
Observe that $y\in B(x,\tfrac{1}{r})$ if and only if $x\in B(y,\tfrac{1}{r})$. Let $m_0>0$ such that $y\in C_{m}$ for all $m\geq m_0$, and let $r_0>0$ such that $B(y,\tfrac{1}{r})\subseteq C_{m_0}$ for all $r\geq r_0$. Then  
for $ r>\max(m_0,r_0)$, it holds  
\begin{equation}\label{Equation:Lemma_Dirac_to_Lebesgue_3}
\int_{C_r} \frac{ \chi_{B(x,1/r)}(y)}{V_{1/r}}f(x)dx =
\int_{C_r} \frac{ \chi_{B(y,1/r)}(x)}{V_{1/r}}f(x)dx= \int_{B(y,1/r)} \frac{ f(x)}{V_{1/r}}dx.
\end{equation}
Since $f$ is continuous at $y$, for a fixed $\epsilon>0$, there exists $\frac{1}{\max(m_0,r_0)}>\eta_\epsilon>0$ such that for all $x\in B(y,\eta_\epsilon)$ it holds $|f(x)-f(y)|<\epsilon$. Thus, for all $\epsilon>0$, there exists $\eta_\epsilon$ such that 	for all $r$ with $\tfrac{1}{r}<\eta_\epsilon$ it holds
	\[  \Big|\Big(\int_{B(y,1/r)} \frac{ f(x)}{V_{1/r}}dx -f(y)\Big)\Big|<\epsilon.\]
This implies that the limit in $r$ of   \eqref{Equation:Lemma_Dirac_to_Lebesgue_3} exists and equals  $f(y)$. This concludes the proof. 
\end{proof}

\begin{lem}\label{Lemma:Kac_Rice_before_random2}
	Let $A\subseteq\R^n$ be a compact set  and $f\colon\R^n\rightarrow\R^n$ a function with continuous first-order partial derivatives in $A$. For $u\in \R^n$, assume 
the solutions to $f(t)=u$ in $A$ are isolated, satisfy $\det\big(J_f(t)\big)\neq 0$, and belong to $A^\circ$. Then
there exists $\lambda>0$ such that for all $\varepsilon<\lambda$ and $v\in B(u,\lambda)$ it holds
	\[\#\big(f^{-1}(v)\cap A\big)= \int_A \frac{\chi_{B(v,\varepsilon)}\big(f(t)\big)}{V_\varepsilon} \, |\det\big(J_f(t)\big)|\, dt.\]
	\end{lem}
\begin{proof}
For $w\in \R^n$, let $S_w$ be the solution set of the equation $f(t)=w$ in $A$, that is $S_w=f^{-1}(w)\cap A$. 
In particular, $S_u$ is finite as the solutions are isolated, and $A$ is compact. 
Furthermore, as $S_u\subseteq A^\circ$ by assumption, there exists $\delta>0$ such that for each $s\in S_u$,  $B(s,\delta) \subseteq A^\circ$ and $B(s,\delta)\cap S_u=\{ s\}$.

As the Jacobian of $f(t)$ does not vanish on the points in $S_u$, by the inverse mapping theorem \cite[Theorem 9.25]{Rudin_Principles}, 
$\delta$ can be chosen such that $f$ is a diffeomorphism from each   $B(s,\delta)$ to $f(B(s,\delta))$.
Choose $\nu>0$ such that $B(u,\nu) \subseteq \bigcap_{s\in S_u} f(B(s,\delta))$.
It follows that for all $v\in B(u,\nu)$, the solutions to $f(t)=v$ also are isolated and the Jacobian of $f(t)$ does not vanish. By choosing $\nu$ smaller if necessary, we further guarantee that all solutions belong to $A^\circ$ as well.  Hence $\#S_u=\#S_v$ and each set $B(s,\delta)$ contains one element of $S_v$.

Let $\lambda=\tfrac{\nu}{3}$ and consider $\varepsilon<\lambda$. 
For $v\in B(u,\lambda)$, 
we have $B(v,\varepsilon)\subseteq B(u,\nu)$ and hence $f$ is a diffeomorphism from each connected component of $f^{-1}(B(v,\varepsilon))$    to $B(v,\varepsilon)$. We denote these connected components by $U_s$, for $s\in S_v$ (which are Borel sets). 
 A change of variables \cite[Theorem 7.26]{Rudin_Real_and_Complex} gives that 
	\begin{align*}
	V_\varepsilon &= \int_{B(v,\varepsilon)}\chi_{B(v,\varepsilon)}(x) dx=\int_{U_s}\chi_{B(v,\varepsilon)}\big(f(t)\big)\, |\det\big(J_f(t)\big)| \, dt.
	\end{align*}
Since  $\chi_{B(v,\varepsilon)}\big(f(t)\big)=0$ if  $t\in A\setminus \cup_{s\in S_v}U_s$, and the union of $U_s$ for $s\in S_v$ is disjoint, by summing over $s\in S_v$ we obtain
\begin{align*}
 \#S_v &= \sum_{s\in S_v} \int_{B(v,\varepsilon)}\frac{ \chi_{B(v,\varepsilon)}(x)}{V_{\varepsilon}} dx 
\\ &= \sum_{s\in S_v}  
\int_{U_s}\frac{\chi_{B(v,\varepsilon)}(f(t))}{V_{\varepsilon}}\, |\det\big(J_f(t)\big)| \, dt +\int_{A\setminus \cup_{s\in S_v}U_s}\frac{\chi_{B(v,\varepsilon)}(f(t))}{V_{\varepsilon}}\, |\det\big(J_f(t)\big)| \, dt \\
	&= \int_A \frac{\chi_{B(v,\varepsilon)}(f(t))}{V_\varepsilon} \, |\det\big(J_f(t)\big)| \, dt.
\end{align*}
Hence, $\#\big(f^{-1}(v)\cap A\big)$ agrees with the integral above, for all $v\in B(u,\lambda)$, and $\varepsilon <\lambda$.
This concludes the proof. 
\end{proof}

\paragraph{\bf Proof of Theorem~\ref{thm:Kac-Rice}.}
We are now ready to proof  Theorem~\ref{thm:Kac-Rice}. 
 As the assumptions of Theorem~\ref{thm:Kac-Rice} hold for a compact box $S\subseteq A$ if they hold for $A$, 
 it is enough to prove Theorem~\ref{thm:Kac-Rice} when $A$ is compact  by Lemma~\ref{Lemma_Kac_Rice_Unbounded_A}. Hence assume $A$ is compact.

	To show that the image of $f_\k(t)$ has positive measure, note that the first $n$ columns of the Jacobian of $f_\k(t)$ with respect to $\k$ form the diagonal matrix with entries $h_i(\bar{\k},t)$, $i=1,\dots,n$, which by assumption do not vanish on $\widetilde{B}\times A$. 

Using \eqref{eq:f_to_g}, the expected value of $\# (f_\k^{-1}(0) \cap A) $ is  given by
		\begin{align*}
	\mathbb{E}(\# (f_\k^{-1}(0) \cap A)) &= \int_{B}   \# (f_\k^{-1}(0) \cap A) \prod_{i=1}^{m} \rho_i(\k_i)    d\k_1\dots d\k_m	 \\
	&= \int_{B}   \# \big(g_{\pi_2(\k)}^{-1}(\pi_1(\k)) \cap A\big) \prod_{i=1}^{m} \rho_i(\k_i)  d\k_1\dots d\k_m =: (\star).
	\end{align*}

By Lemma~\ref{lem:refined}(ii), outside a measure zero set of $B$, the equation $f_\k(t)=0$ has 
a finite number of solutions. 
As $f_\k(t)$ is polynomial, there is an upper bound $M$ on the number of complex solutions  depending only on the exponents of the monomials, and not on the coefficients (this follows for instance from Bernstein-Kushnirenko theorem  on the number of solutions in the torus $(\mathbb{C}^*)^n$ \cite[Theorem 5.4]{Cox}).  Hence the integrand in $(\star)$, which is non-negative, is bounded above by  $M$ and the integral $(\star)$   is finite.  

Let $\rho^{(1)}$, $\rho^{(2)}$ be the densities on $\pi_1(B)$ and $\pi_2(B)$, respectively. 
 By Tonelli's theorem,  the integral over $B$ can be found iteratively over $\pi_2(B)=\widetilde{B}$ and $\pi_1(B)$ with variables $y,z$ respectively:
		\begin{align*}
	 (\star) &= \int_{\widetilde{B}}\left( \int_{\pi_1(B)}  \# \big(g_{y}^{-1}(z) \cap A\big)	 \rho^{(1)}(z) dz \right) \rho^{(2)}(y) dy.
	 \end{align*}

By Lemma~\ref{lem:refined}(i), 
for every $y$ outside a measure zero set $\widetilde{P}'$, 
 there exists a measure zero set $A(y)$ such that $\rho_i$ is continuous at $g_i(y,t)$ for all $t\notin A(y)$ and for all $i=1,\dots,n$. 
 Let $\widetilde{P} $ be the set in Lemma~\ref{lem:refined}(ii), consider the (relative open) set $\widetilde{B}':= \widetilde{B}\setminus (\widetilde{P} \cup \widetilde{P}')$, and fix $y\in \widetilde{B}'$.
We focus on the inner integral: 
\begin{align*}
(\star\star) &:=\int_{\pi_1(B)}  \# \big(g_{y}^{-1}(z) \cap A\big)	 \rho^{(1)}(z) dz.
\end{align*}
As the denominator of $g_{y}$ does not vanish on $\widetilde{B}\times A$, $g_{y}$  has continuous first-order partial derivatives in $t\in A$ for all $y\in \widetilde{B}$. 
By Lemma~\ref{lem:refined}(ii), there exists a measure zero  set $Q_y$ such that for $z\in \pi_1(B)\setminus Q_y$, the solutions to $g_y(t)=z$ are isolated, belong to  $A^\circ$, and the Jacobian does not vanish. Lemma~\ref{Lemma:Kac_Rice_before_random2} applies to $g_{y}$, with $u=z\in \pi_1(B)\setminus Q_y$. 
Hence, for every such pair $(z,y)$ and $\varepsilon>0$ small enough, we have 
\begin{equation}\label{eq:niceeq}
\#(g_{y}^{-1}(z)\cap A)
 = \int_A\frac{\chi_{B(z,\varepsilon)}(g_{y}(t))}{V_\varepsilon}|\det(J_{g_{y}}(t))| \, dt.
 \end{equation}

Let $B_1' :=\pi_1(B)\setminus (Q_y\cup \partial \pi_1(B))$, which is an open set, 
and consider the integral $(\star\star)$ over $B_1'$ instead. As $Q_y\cup \partial \pi_1(B)$ has measure zero, the value of the integral is the same.  For every $r\in\N$, let $D_r\subseteq B_1'$ be the set of all points $z\in B_1'$ for which 
\eqref{eq:niceeq} holds for $\varepsilon\leq \tfrac{1}{r}$, and let $C_r$ be its interior. 
Clearly $C_1\subseteq C_2\subseteq \dots$ defines an increasing sequence of open sets. Their  union is 
\[\lim_{r\rightarrow +\infty}C_r=B_1' \]
by Lemma~\ref{Lemma:Kac_Rice_before_random2},  as if $z\in B_1'$, there exists $r_0>0$ and $\delta>0$ such that $B(z,\delta)\subseteq D_{r_0}$, hence $z\in C_{r_0}$. 
The above discussion, \eqref{eq:niceeq} and Lebesgue's monotone convergence theorem \cite[Theorem 1.26]{Rudin_Real_and_Complex} gives that
\begin{align*}
(\star\star)  &=  \lim_{r\rightarrow +\infty}\int_{C_r} \# (g_{y}^{-1}(z) \cap A)\, \rho^{(1)}(z) dz 
\\ &= \lim_{r\rightarrow +\infty}\int_{C_r}\left( \int_A\frac{\chi_{B(z,1/r)}\big(g_{y}(t)\big)}{V_{1/r}}\, |\det\big(J_{g_{y}}(t)\big)| dt \, \right)  \rho^{(1)}(z) dz.  \end{align*}
We now claim that the following equalities, derived from interchanging limits, hold:
\begin{align}
(\star\star)  &
= \lim_{r\rightarrow +\infty} \int_A \left(\int_{C_r}
\frac{\chi_{B(z,1/r)}\big(g_{y}(t)\big)}{V_{1/r}}\, \rho^{(1)}(z)dz
 \right) |\det\big(J_{g_{y}}(t)\big)|  \, dt  \label{eq:eq1}
\\ &=
 \int_A \lim_{r\rightarrow +\infty} \left(\int_{C_r} \frac{\chi_{B(z,1/r)}\big(g_{y}(t)\big)}{V_{1/r}}\, \rho^{(1)}(z)dz
 \right) |\det\big(J_{g_{y}}(t)\big)|  dt \label{eq:eq2}
\\ &=  \int_A    \rho^{(1)}(g_{y}(t)) \,  |\det\big(J_{g_{y}}(t)\big)|  \, dt.
\label{eq:eq3}
\end{align}

Let us show that \eqref{eq:eq1}-\eqref{eq:eq3} hold. Recall  $y\in \widetilde{B}'$ is fixed. 
The interchange of limits in \eqref{eq:eq1} follows again from Tonelli's theorem \cite[Theorem 22.7]{Aliprantis-Burkinshaw-1990}, as the integrand is a non-negative measurable function.
For \eqref{eq:eq2}, we need to show that we can interchange the limit in $r$ and the integral in $A$. To this end, 
we   show that Lebesgue's dominated convergence theorem \cite[Theorem 1.34]{Rudin_Real_and_Complex} applies. Consider the sequence 
\[\alpha_r := \int_{C_r} \frac{\chi_{B(z,1/r)}\big(g_{y}(t)\big)}{V_{1/r}}\, \rho^{(1)}(z)dz. \]
Let $0<M_1<+\infty$ be an upper bound of $\rho^{(1)}$. We have 
\begin{align*}
\alpha_r &  \leq  
\int_{\pi_1(B)} \frac{\chi_{B(z,1/r)}\big(g_{y}(t)\big)}{V_{1/r}}\, \rho^{(1)}(z)dz  =\int_{\pi_1(B)} \frac{\chi_{B(g_{y}(t),1/r)}(z)}{V_{1/r}}\, \rho^{(1)}(z)dz 
\\
&= \int_{B(g_{y}(t),1/r)\cap \pi_1(B)} \frac{1}{V_{1/r}}\, \rho^{(1)}(z)dz 
 \leq M_1.
 \end{align*}
 As $|\det\big(J_{g_{y}}(t)\big)|$ is integrable, the sequence $\alpha_r \, |\det\big(J_{g_{y}}(t)\big)|$ is dominated by an integrable function. 
This gives  \eqref{eq:eq2}, as long as the limit exists, but this   follows from the proof of \eqref{eq:eq3}.

For  \eqref{eq:eq3}, we apply Lemma~\ref{Lemma:Dirac_to_Lebesgue3} 
with $f=\rho^{(1)}$, point $g_y(t)$ with fixed  $t$ and set $T=B_1'$. 
Note that it is enough to prove that  \eqref{eq:eq3} holds for the integral over $A\setminus A(y)$ instead of over $A$, as $A(y)$ has measure zero. To this end, it is enough to verify that the hypotheses of Lemma~\ref{Lemma:Dirac_to_Lebesgue3}  hold.  First, the choices made above imply that  for  $t\notin A(y)$,  $\rho^{(1)}$ is continuous at $g_y(t)$, and $g_y(t)\notin \partial B_1'\subseteq \partial \pi_1(B) \cup Q_y$. The last condition 
follows from Lemma~\ref{Lemma:Kac_Rice_before_random2}.

\medskip
Finally, using the  expression found for $(\star\star)$, the definition of $\bar{\rho}$ in the statement, and that $\widetilde{P}\cup \widetilde{P}'$ has measure zero, we obtain
\begin{align*}
(\star) &= \int_{\widetilde{B}} \left(  \int_A  \,  |\det\big(J_{g_{y}}(t)\big)| \,  \rho^{(1)}(g_{y}(t))  \, dt \right)\rho^{(2)}(y) dy .
\end{align*}
All that is left is to justify that the integrals in $A$ and $\widetilde{B}$ can be interchanged, but, again, this follows from Tonelli's theorem, as the integrand is non-negative and measurable. 
This concludes the proof of Theorem~\ref{thm:Kac-Rice}. \hfill\qed

  \section{Conclusion}
For a type of parametric polynomial systems of equations, this work has proposed a numerical approach to partition the parameter region into the regions where the number of solutions of the system is constant.
The idea builds on classical  Kac-Rice formulas for the number of zeroes of random functions. In particular, we derive a Kac-Rice formula for the average number of solutions that the system of polynomial equations has in a  box $A$, when parameters follow given probability distributions. After partitioning the parameter region into boxes, and applying the Kac-Rice formula in each box, we obtain a coarse approximation of the desired regions.

The numerical aspect of the approach resides in the computation of the integral of the Kac-Rice formula. We propose to use Monte-Carlo integration and exploit the fact that we compute the integral over $A$ of the expected value of a function in some of the parameters. For the integral over $A$, there is no obvious generic choice of probability distribution to be used with Monte-Carlo integration. 

We have demonstrated through several examples how the approach can be successfully applied. As in our examples the box $A$ is bounded, the choice of uniform distribution on $A$  turned out to be suitable.
A detailed analysis of the complexity of our approach is out of the scope of this work, but the examples illustrate how the method can handle relatively large number of parameters and variables. 
The main limitations are posed by the sample size required to compute,  with the desired precision, the Kac-Rice integral using Monte-Carlo integration, and, additionally, by the number of integrals to be computed.

\bigskip
\textbf{Acknowledgements. }  
The authors acknowledge funding from the Independent Research Fund of Denmark. 
This work was initiated while A. S. visited MPI for the Mathematical Sciences in Leipzig in the summer of 2017. In particular A. S. learned about the Kac-Rice formula in the \emph{Reading Group on Real Algebraic Geometry} that took place at MPI in June 2017. We thank Paul Breiding for clarifications on the Kac-Rice formula and Bernd Sturmfels for discussions on algebraic approaches to determine parameter regions of multistationarity.
	We thank Jimmy Olsson for discussions on Monte Carlo methods for numerical integration, 
Carsten Wiuf for key and fruitful discussions on preliminary drafts of this manuscript, and Matthew England for discussions on algorithms in semi-algebraic geometry.

	{\small

%\bibliographystyle{plain}
%\bibliography{binomiality}
}

\end{document}